\documentclass[11pt,reqno]{amsart}
\usepackage{amsmath}
\usepackage{amssymb}
\usepackage{amsfonts}

\def\squarebox#1{\hbox to #1{\hfill\vbox to #1{\vfill}}}

\newcommand{\Z}{{\mathbb Z}}

\newcommand{\R}{{\mathbb R}}
\newcommand{\C}{{\mathbb C}}
\newcommand{\N}{{\mathbb N}}

\theoremstyle{plain}

\newtheorem{thm}{Theorem}

\newtheorem{lem}{Lemma}
  
\newtheorem{prop}{Proposition}

\newtheorem{deff}{Definition}

\newtheorem{rem}{Remark}

\begin{document}

\def\sn{{\bf S}^{n-1}}
\def\ts{\tilde{\sigma}}
\def\ss{{\mathcal S}}
\def\aa{{\mathcal A}}
\def\cc{{\mathcal C}}
\def\tpp{\widetilde{{\mathcal P}}}
\def\iu{\underline{i}}
\def\lc{{\mathcal L}}
\def\pxi{\phi + (\xi + i u)\psi_p}
\def\pp{\mathcal P}
\def\rr{\mathcal R}
\def\hU{\widehat{U}}
\def\mt{\Lambda}
\def\hLa{\widehat{\mt}}
\def\ep{\epsilon}
\def\tPhi{\widetilde{\Phi}}
\def\hR{\widehat{R}}
\def\oV{\overline{V}}
\def\uu{\mathcal U}
\def\tz{\tilde{z}}
\def\hz{\hat{z}}
\def\hd{\hat{\delta}}
\def\ty{\tilde{y}}
\def\hs{\hat{s}}
\def\hc{\hat{\cc}}
\def\hC{\widehat{C}}
\def\hh{\mathcal H}
\def\tf{\tilde{f}}
\def\of{\overline{f}}
\def\trr{\tilde{r}}
\def\tr{\tilde{r}}
\def\tts{\tilde{\sigma}}
\def\tVl{\widetilde{V}^{(\ell)}}
\def\tVj{\widetilde{V}^{(j)}}
\def\tVo{\widetilde{V}^{(1)}}
\def\tVj{\widetilde{V}^{(j)}}
\def\tPsi{\tilde{\Psi}}
\def\tp{\tilde{p}}
\def\tVjo{\widetilde{V}^{(j_0)}}
\def\tvj{\tilde{v}^{(j)}}
\def\tVjj{\widetilde{V}^{(j+1)}}
\def\tvl{\tilde{v}^{(\ell)}}
\def\tVt{\widetilde{V}^{(2)}}
\def\Lo{\; \stackrel{\circ}{L}}
\def\tg{\tilde{g}}

\def\ii{{\imath }}
\def\jj{{\jmath }}
\vspace*{0,8cm}
\def\saa{\Sigma_A^+}
\def\sAA{\Sigma_{\aa}^+}
\def\Lip{\mbox{\rm Lip}}
\def\clip{C^{\mbox{\footnotesize \rm Lip}}}
\def\lip{\mbox{{\footnotesize\rm Lip}}}
\def\Vol{\mbox{\rm Vol}}

\def\ccm{\cc^{(m)}}
\def\oo{\mbox{\rm O}}
\def\ooo{\oo^{(1)}}
\def\oot{\oo^{(2)}}
\def\ooj{\oo^{(j)}}
\def\fo{f^{(0)}}
\def\ft{f^{(2)}}
\def\fj{f^{(j)}}
\def\wo{w^{(1)}}
\def\wt{w^{(2)}}
\def\wj{w^{(j)}}
\def\Vo{V^{(1)}}
\def\Vt{V^{(2)}}
\def\Vj{V^{(j)}}

\def\Uo{U^{(1)}}
\def\Ut{U^{(2)}}
\def\Ul{U^{(\ell)}}
\def\Uj{U^{(j)}}
\def\wl{w^{(\ell)}}
\def\Vl{V^{(\ell)}}
\def\Ujj{U^{(j+1)}}
\def\wjj{w^{(j+1)}}
\def\Vjj{V^{(j+1)}}
\def\Ujo{U^{(j_0)}}
\def\wjo{w^{(j_0)}}
\def\Vjo{V^{(j_0)}}
\def\vj{v^{(j)}}
\def\vl{v^{(\ell)}}

\def\gl{\gamma_\ell}
\def\id{\mbox{\rm id}}
\def\piU{\pi^{(U)}}
\def\bs{\bigskip}
\def\ms{\medskip}
\def\Int{\mbox{\rm Int}}
\def\diam{\mbox{\rm diam}}
\def\di{\displaystyle}
\def\dist{\mbox{\rm dist}}
\def\ff{{\cal F}}
\def\i{{\bf i}}
\def\pr{\mbox{\rm pr}}
\def\co{\; \stackrel{\circ}{C}}
\def\la{\langle}
\def\ra{\rangle}
\def\supp{\mbox{\rm supp}}
\def\Arg{\mbox{\rm Arg}}
\def\Int{\mbox{\rm Int}}
\def\II{{\mathcal I}}
\def\tll{\widetilde{\lc}}
\def\e{\emptyset}
\def\endofproof{{\rule{6pt}{6pt}}}
\def\con{\mbox{\rm const }}
\def\Box{\spadesuit}
\def\be{\begin{equation}}
\def\ee{\end{equation}}
\def\beqn{\begin{eqnarray*}}
\def\eeqn{\end{eqnarray*}}
\def\MM{{\mathcal M}}
\def\tmu {\tilde{\mu}}
\def\Pr{\mbox{\rm Pr}}
\def\Prf{\mbox{\footnotesize\rm Pr}}
\def\htau{\hat{\tau}}
\def\btau{\overline{\tau}}
\def\hr{\hat{r}}
\def\tF{\widetilde{F}}
\def\tG{\widetilde{G}}
\def\trho{\tilde{\rho}}

\title[Spectral estimates and sharp large deviations]{Spectral estimates for Ruelle operators\\ with two parameters and sharp large deviations}
\author[V. Petkov]{Vesselin Petkov}
\address{Universit\'e de Bordeaux, Institut de Math\'ematiques de Bordeaux, 351,
Cours de la Lib\'eration, 
33405  Talence, France}
\email{petkov@math.u-bordeaux.fr}
\author[L. Stoyanov]{Luchezar Stoyanov}
\address{University of Western Australia, Department of Mathematics 
and Statistics, Perth, WA 6009, Australia}
\email{luchezar.stoyanov@uwa.edu.au}

\def\ts{\tilde{\sigma}}
\def\ss{{\mathcal S}}
\def\aa{{\mathcal A}}
\def\R{{\mathbb R}}
\def\C{{\mathbb C}}
\def\iu{\underline{i}}
\vspace*{0,8cm}
\def\saa{\Sigma_A^+}
\def\sAA{\Sigma_{\aa}^+}
\def\lc{{\mathcal L}}
\def\pxi{\phi + (\xi + i u)\psi_p}
\def\sa{\Sigma_A}
\def\ssa{\Sigma_A^+}
\def\ssn{\sum_{\sigma^n x = x}}
\def\lc{{\mathcal L}}
\def\lo{{\mathcal O}}
\def\oo{{\mathcal O}}
\def\mt{\Lambda}
\def\ep{\epsilon}
\def\ms{\medskip}
\def\bs{\bigskip}
\def\diam{\mbox{\rm diam}}
\def\rr{\mathcal R}
\def\pp{\mathcal P}
\def\hU{\widehat{U}}
\def\hz{\hat{z}}
\def\tz{\tilde{z}}
\def\ty{\tilde{y}}
\def\i{{\bf i}}
\def\pp{{\mathcal P}}
\def\ff{{\mathcal F}_{\theta}}
\def\tG{\tilde{G}}
\def\hw{\hat{w}}
\def\ep{\epsilon}
\def\rt{R^{\tau}}
\maketitle

\def\tm{\tilde{m}}
\def\tj{\tilde{j}}
\def\dd{{\mathcal D}}
\def\piU{\pi^{(U)}}
\def\hrho{\hat{\rho}}
\def\hdd{\widehat{\dd}}
\def\Xijl{X^{(\ell)}_{i,j}}
\def\hXijl{\widehat{X}^{(\ell)}_{i,j}}
\def\omijl{\omega^{\ell}_{i,j}}
\def\vl{v^{(\ell)}}
\def\kk{{\mathcal K}}
\def\Lip{\mbox{\rm Lip}}
\def\ff{{\mathcal F}}
\def\diamte{\diam_\theta}
\def\J{{\sf J}}
\def\nn{{\mathcal N}}
\def\ma{{\mathcal M}_a}
\def\mac{{\mathcal M}_{ac}}
\def\dte{D_\theta}
\def\fa{f^{(a)}}
\def\lab{{\mathcal L}_{ab}}
\def\dl{d^{(\ell)}}
\def\tu{\tilde{u}}
\def\hu{\hat{u}}
\def\hz{\hat{z}}
\def\Gl{\Gamma^{(\ell)}}
\def\lambdam{\lambda^{(m)}}
\def\thetam{\theta^{(m)}}
\def\tR{\widetilde{R}}
\def\htheta{\hat{\theta}}
\def\hW{\widehat{W}}
\def\hnu{\hat{\nu}}
\def\labw{\lc_{abw}}
\def\labz{\lc_{abz}}
\def\tcc{\tilde{\cc}}
\def\tdd{\widetilde{\dd}}
\def\hcc{\widehat{\cc}}
\def\lambdam{\lambda^{(m)}}
\def\thetam{\theta^{(m)}}
\def\Lipf{\mbox{\rm\footnotesize Lip}}
\def\yl{y^{(\ell)}}
\def\tS{\widetilde{S}}
\def\E{\mathcal E}
\def\tp{\tilde{p}}
\def\tq{\tilde{q}}
\def\f0{f^{(0)}}
\def\ttau{\tilde{\tau}}

\def\fatc{f_{atc}}
\def\gt{g_t}
\def\Gt{G_t}
\def\tft{\tilde{f}_t}
\def\ft{f_t}
\def\fat{f_{at}}
\def\MM{{\mathcal M}}
\def\mm{{\mathcal M}}
\def\fa{f_a}
\def\fb{f+b}
\def\f0{f_0}
\def\fab{f_{ab}}
\def\mab{\mm_{ab}}
\def\ma{\mm_{a}}
\def\lab{\lc_{ab}}
\def\psib{\psi^{(b)}}
\def\fac{f_{ac}}
\def\mac{\mm_{ac}}
\def\mat{\mm_{at}}
\def\matc{\mm_{atc}}
\def\hatc{h_{atc}}
\def\lambdaatc{\lambda_{atc}}
\def\tmatc{\widetilde{\mm}_{atc}}
\def\tmat{\widetilde{\mm}_{at}}
\def\lac{\lc_{ac}}
\def\labc{\lc_{abc}}
\def\labt{\lc_{abt}}
\def\labtz{\lc_{abtz}}
\def\tlabtz{\tll_{abtz}}
\def\labz{\lc_{abz}}
\def\psic{\psi^{(c)}}

\def\hclip{\hat{C}^{\lip}}
\def\hC{\hat{C}}

\def\PPsi{{\sf \Psi}}
\def\hk{\hat{k}}
\def\hgamma{\hat{\gamma}}

\def\Intu{\mbox{\rm Int}^u}
\def\hZ{\widehat{Z}}
\def\xijl{X^{(\ell)}_{i,j}}
\def\eijl{\omega^{(\ell)}_{i,j}}
\def\hxijl{\widehat{X}^{(\ell)}_{i,j}}

\def\gao{\gamma^{(1)}}
\def\gat{\gamma^{(2)}}
\def\diamtef{{\footnotesize\mbox{\rm  diam}_\theta}}
\def\Con{\mbox{\rm Const}}
\def\Conf{\mbox{\footnotesize\rm Const}}

\def\gej{\chi^{(j)}_\mu}
\def\ge{\chi_\epsilon}
\def\geo{\chi^{(1)}_\mu}
\def\get{\chi^{(2)}_\mu}
\def\gei{\chi^{(i)}_{\mu}}
\def\gee{\chi_{\mu}}
\def\gett{\chi^{(2)}_{\mu}}
\def\geol{\chi^{(1)}_{\ell}}
\def\getl{\chi^{(2)}_{\ell}}
\def\geil{\chi^{(i)}_{\ell}}
\def\gee{\chi_{\ell}}

\def\hmu{\hat{\mu}}
\def\tgt{\tilde{g}_t}
\def\tGt{\widetilde{G}_t}
\def\tPsit{\widetilde{\Psi}_t}
\def\stt{\sigma_{\tau}^t}
\begin{abstract}
We obtain spectral estimates for the iterations of Ruelle operator $L_{f + (a + \i b)\tau + (c + \i d) g}$ with two complex parameters and H\"{o}lder functions $f,\: g$ generalizing the case $\Pr(f) =0$ studied in \cite{PeS2}. As an application we prove a sharp large deviation theorem concerning exponentially shrinking intervals which  improves the result in \cite{PeS1}.

\end{abstract}

\section{Introduction}
\renewcommand{\theequation}{\arabic{section}.\arabic{equation}}
\setcounter{equation}{0}

Let $M$ be a $C^2$ complete Riemannian manifold,  let $\varphi_t : M \longrightarrow M$ ($t\in \R$) be a $C^2$  flow on $M$ 
and let $\varphi_t: M \longrightarrow M$  be a $C^2$ weak mixing Axiom A flow (\cite{KH}, \cite{PP}). 
Let $\mt$ be a {\it basic set} for $\phi_t$, that is, $\mt$ is a compact
locally maximal invariant subset of $M$ and  $\varphi_t$ is hyperbolic and transitive on $\mt$. 

As in \cite{PeS2}, we will use a symbolic coding of the flow on $\Lambda$ provided by a
 a fixed Markov family $\{ R_i\}_{i=1}^k$. More precisely, we consider a Markov family of 
{\it pseudo-rectangles}  $R_i = [U_i  , S_i ] =  \{ [x,y] : x\in U_i, y\in S_i\}$ (see section 2 for more details). 
Denote by  $\pp: R = \cup_{i=1}^{k} R_i \longrightarrow R$ the related Poincar\'e map, by
$\tau(x) > 0$ the first return time function on $R$, and by  $\sigma : U = \cup_{i=1}^k U_i \longrightarrow U$ the {\it shift map} 
given by $\sigma  = \piU \circ \pp$, where $\piU : R \longrightarrow U$ is the {\it projection} along stable leaves. 
The flow $\varphi_t$ on $\mt$ is naturally related to the suspension flow  $\sigma_t^{\tau}$ on the suspension space
$R^{\tau}$ (see section 2 for details). There exists a natural semi-conjugacy projection $\pi(x, t): R^{\tau} \longrightarrow \Lambda$ 
which is one-to-one on a residual set (see \cite{B2}). For $x \in R$ set
$$\tau^n(x): = \tau(x) + \tau(\sigma(x)) +...+ \tau(\sigma^{n-1}(x)).$$

Given H\"older continuous functions $F, G : \Lambda \longrightarrow \R$, define $f,g : R \longrightarrow \R$ by
$$f(x) = \int_0^{\tau(x)} F(\pi(x, t)) dt \quad , \quad g(x) = \int_0^{\tau(x)} G(\pi(x, t)) dt .$$
The main object of study in this paper are the Ruelle transfer operators of the form
$$L_{f - s\tau + zg} v(x) = \sum_{\sigma y = x} e^{f(y) - s \tau(y) + zg(y)} v(y) \quad, \quad s, z \in \C \quad , x \in U ,$$
depending on two complex parameters $s$ and $z$. Under certain assumptions, strong spectral estimates for such operators
have been established in \cite{PeS2} and some significant applications to the study of zeta functions depending on two
complex parameters have been made. 
 We denote by $m_H$ the equilibrium state corresponding to $H$ in $R^{\tau}$ and by $\mu_k$ the equilibrium state corresponding to $k$ in $R$. 
 More precisely,
$$\Pr_{\sigma}(k) = h(\sigma, \mu_k) + \int_{R} k d\mu_k = \sup_{\mu \in {\mathcal M}_{\sigma}} \Bigl\{ h(\sigma, \mu) + \int_R k d\mu\Bigr\},$$
$$\Pr_{\sigma_{\tau}}(H) = h(\sigma_{\tau}^t, m_H) + \int_{R^{\tau}} H dm_H = \sup_{m \in {\mathcal M}_{\sigma_{\tau}}}\Bigl\{ h(\sigma_{\tau}^t,m) + \int_{R^{\tau}} H dm \Bigr\},$$
where $h(\sigma, \mu)$ is the metric entropy of $\sigma$ with respect to $\mu$ and $h(\sigma_{\tau}^t, m)$ is the metric entropy of the suspended flow 
$\sigma_{\tau}^t$ with respect to $m$. Let $P = \Pr_{\sigma}(f)$.

Let $\| h\|_0$ denote the {\it standard $\sup$ norm} of $h$ on $U$.  For  $|b| \geq 1$, and $\beta > 0$, as in \cite{D}, 
define the norm  $\|h\|_{\beta,b} = \|h\|_\infty + \frac{|h|_\beta}{|b|}$ on the space $C^\beta(U)$ of $\beta$-H\"older
functions on $U$.

Our first aim in this paper is to prove the following theorem.

\ms

\begin{thm} 
Let $\varphi_t : M \longrightarrow M$ satisfy the Standing Assumptions (see Sect. $4$) over the basic set $\mt$, and let $0 < \beta < \alpha$.
Let  $\rr = \{R_i\}_{i=1}^k$ be a Markov family for $\varphi_t$ over $\mt$ as in section 2. Then 
for any real-valued functions  $f,g \in C^\alpha(\hU)$ and any constants $\epsilon > 0$ and $B > 0$ 
there exist constants $0 < \rho < 1$, $a_0 > 0$, $b_0 \geq 1$ and  $C = C(B, \epsilon)> 0$ 
such that if $a,c\in \R$ satisfy $|a|, |c| \leq a_0$ then
\begin{equation} \label{eq:1.1}
\|L_{f -(a+ \i b)\tau + (c+\i w) g}^m h \|_{\beta,b}  \leq C \,e ^{Pm}\, \rho^m \, |b|^{\ep}\, \| h\|_{\beta,b}
\end{equation}
for all $h \in C^\beta(U)$, all integers $m \geq 1$ and all $b, w\in \R$ with  $|b| \geq b_0$ and $|w| \leq B \, |b|$.

\end{thm}

\ms

In Theorem 5.1 in \cite{PeS2} the above estimate has been proved in the case $P = 0$ assuming $|w| \leq B \, |b|^\nu$ for some constant $\nu \in (0,1)$.
The present results is significantly stronger. See also Remark 1 below.

In the proof of Theorem 1 we will use some arguments from the proof of Theorem 5.1 in \cite{PeS2} with necessary modifications.

\begin{rem}
Notice that in Theorem $1$ above we do not assume that pressure $P$ of $f$ is zero, unlike what has been done in
previous papers. This contributes the term $e^{Pm}$ in the right-hand-side of $(1.1)$ which is significant especially in
the case $P < 0$ which occurs in the applications concerning large deviations (see Section $3$). In previous papers the 
authors consider the case $P = 0$ and remark that the general case follows 
from this. However a more careful argument shows that an estimate of the form $(1.1)$ does not follow immediately from
a similar estimate\footnote{Indeed, assume we have proved $(1.1)$ in the case $P = 0$, and then deal with the general case using the standard approach.
Given $a,b$ as in the theorem and  $h\in \ff_\theta (\hU)$, we have
\begin{eqnarray} \label{eq:1.2}
(L_{f -(a+ \i b)\tau}^m h) (u) 
 =  \sum_{\sigma^m(v) = u} e^{(f- (a+\i b)\tau)^m(v)} h(v) =  (L_{f -(P+a+ \i b)\tau}^m (e^{P\tau^m(v)}\, h)) (u) .
\end{eqnarray}
We can now apply $(1.1)$ in the case $P = 0$ replacing $h$ by $e^{P\tau^m(v)}\, h$. 
Since $0 < c \leq \tau(u) \leq c_1$ for some constants $c$ and $c_1$, assuming e.g.
$P \leq 0$ (the other case is similar), we get
$\| e^{P\tau^m}\, h\|_\infty \leq e^{mPc} \|h\|_\infty ,$
and
$| e^{P\tau^m}\, h|_\theta \leq e^{m P c} \,  |h|_\theta + |e^{P\tau^m}|_\theta \, \|h\|_\infty  .$
Given $u,v\in U_i$ we have (assuming e.g. $e^{P \tau^m(u)} > e^{P\tau^m(v)}$)
\begin{eqnarray*}
\left| e^{P \tau^m(u)} - e^{P\tau^m(v)}\right|
 \leq  e^{P \tau^m(u)} \; |P \tau^m(u) - P\tau^m(v) |
\leq |P|\, e^{mPc}\, \sum_{j=0}^{m-1} \frac{D_\theta(u,v)}{\theta^j} \leq  \Con \frac{e^{mPc}}{\theta^m}  \, D_\theta(u,v) .
\end{eqnarray*}
Taking $\theta$ closer to $1$ and replacing $c$ by some $c_0< c$, we get $|e^{P\tau^m}|_\theta \leq \Con \; e^{mPc_0}$.
This implies
$$\| e^{P\tau^m}\, h\|_{\theta, b} \leq e^{mPc} \|h\|_\infty + \frac{1}{|b|} \left(e^{m P c} \,  |h|_\theta 
+ \Con \; e^{mPc_0} \; \|h\|_\infty\right) \leq \Con\; e^{mPc_0}\; \|h\|_{\theta,b} .$$
Combining the latter with (1.2) gives
$\|L_{f -(a+ \i b)\tau}^m h \|_{\theta,b} 
 \leq  C \;\rho^m \; \| e^{P\tau^m(v)}\, h \|_{\theta,b} \leq C \; e^{mPc_0}\; \rho^m \;  \|h\|_{\theta,b} .$
As one can see this estimate is a bit worse than (1.1), since $c > c_0 > 0$ (and also $c_1$ in the case $P >0$)
can be rather small constants.} with $P = 0$. 
\end{rem}

\ms
\begin{rem} In the proof of Theorem $2$ in Section $3$ we apply Theorem $1$ with $b = C w$ for some  constant $C > 0$; then $|w| = \frac{1}{C} |b|$.
The relevant part of Theorem $5.1$ in \cite{PeS2} assumes $|w| \leq B |b|^{\nu}$ for some $\nu \in (0, 1)$ and this is clearly not sufficient for the proof
of Theorem $2$ below.
\end{rem}

Let $G$ be a H\"older function on $\Lambda$ such that $G > 0$ everywhere on $\Lambda$. Consider a number
$$0 < a = \int_{R^{\tau}} G dm_{F + \xi(a) G} \in \Bigl\{\int_{R^{\tau}} G dm_{F + t G},\:t \in \R\Bigr\},$$
where $\xi(a)$ is determined by the equation
$$\frac{d \Pr_{\sigma_{\tau}}(F +t G)}{dt}\Bigr \vert_{t = \xi(a)} - a = 0.$$
Let $0 < \rho < 1$ be the constant from Theorem 1, and let $0 < \alpha_0 = - \frac{\log \rho}{2}$. Fix
an arbitrary $0 < \delta \leq \alpha_0$ and consider the sequence $\{\delta_n\}_{n \in \N}$, where
$$\delta_n = e^{-\delta n} .$$   
Set $g_a = g- \tau a$. Then
$$g_a^n(x) = g^{n}(x) - \tau^n(x) a =  \int_0^{\tau^n(x)}G(\pi(t, x))dt - \tau^{n}(x) a.$$
Clearly the property  
$$\frac{g^n(x)}{\tau^n(x)} - a \in \Bigl(-\frac{\delta_n}{\tau^n(x)}, \frac{\delta_n}{\tau^n(x)}\Bigr),$$
is equivalent to 
$$g^n(x) - \tau^n(x) a \in (-\delta_n,\delta_n).$$
On the other hand, since $c n \leq\tau^n(x)\leq c_1 n,\: \forall x \in R, \forall n \in \N$ with some constants $0 < c \leq c_1$ for every $x$,  the
interval $\Bigl(-\frac{\delta_n}{\tau^n(x)}, \frac{\delta_n}{\tau^n(x)}\Bigr)$ is exponentially shrinking to 0 as $n \to \infty.$ Let $\mu =\mu_f$ 
be the equilibrium state of $f$. 

Our second problem concerns the analysis of the asymptotic of
$$\mu\{x: g^n(x) - \tau^n(x)a \in (-\delta_n, \delta_n)\},\:n \to \infty$$
and for $a \neq \int_{R^{\tau}} G dm_F$ we obtain a large deviation result. On the other hand, as in the previous paper
\cite{PeS1}, we examine the measure of points $x \in R$ for which the difference $\frac{g^n(x)}{\tau^n(x)} - a$ stays in an exponentially shrinking interval. 
Next we state two definitions from  \cite{La} and \cite{W}.

\begin{deff} Two functions $f_1, f_2$ are called $\sigma-$independent if whenever there are constants $t_1, t_2 \in \R$ such that
$t_1 f_1 + t_2 f_2$ is cohomologous to a function in $C(R : 2 \pi \Z)$ , we have $t_1 = t_2 = 0.$
\end{deff}

For a function $G \in C^{\beta}(R^{\tau}: \R)$ consider the skew product flow $S_t^G$ on ${\mathbb S}^1 \times R^{\tau}$ defined by
$$S_t^G ( e^{2 \pi \i \alpha}, y) = \Bigl(e^{2 \pi \i (\alpha + \int_0^t G (\varphi_{\tau}^sy) ds)}, \stt(y)\Bigr).$$

\begin{deff} Let $G \in C^{\beta}(R^{\tau}: \R))$. Then $G$ and $\stt$ are {\it flow independent} if 
the following condition is satisfied: if $t_0, t_1 \in \R$ are constants such that the skew product flow $S_t^H$ with $H = t_0 + t_1 G$ 
is not topologically ergodic, then $t_0 = t_1 = 0.$
\end{deff}

Notice that if $G$ and $\stt$ are flow independent, the flow $\stt$ is topologically weak mixing and the function $G$ is not 
cohomologous to a constant function. This implies that the set
$$\Bigl\{\int_{R^{\tau}} G dm_{F + t G},\:t \in \R\Bigr\}$$
has a non-empty interior and setting $\beta(t) = \Pr_{\sigma_{\tau}}(F + t G),$ one has
$$\beta''(t) = \frac{d^2 \Pr_{\sigma_{\tau}}(F + t G)}{dt^2} = \sigma^2_{m_{F + tG}}(G)$$
with 
$$\sigma_m^2(G) = \lim_{T \to \infty}\frac{1}{T}\Bigl\{\int_{R^{\tau}}\int_0^T G(\stt(y)) dtdm - T \int_{R^{\tau}} G dm\Bigr\}^2 < \infty.$$
Moreover, $\beta'(\xi(a)) = a$ and $\xi(a)$ is differentiable with  $\xi'(a) = \frac{1}{\beta''(\xi(a))} > 0$. Without loss of generality by adding a constant, we may assume 
that $\Pr_{\sigma_{\tau}}(F) = 0$. Then $m_F$ and $\Pr_{\sigma_{\tau}}(F + t G)$ don't change and $\Pr_{\sigma}( f -\Pr_{\sigma_{\tau}}(F) \tau)= 0$ yields $\Pr_{\sigma}(f) = 0$. Introduce the rate function
$$\gamma(a)=: \Pr_{\sigma_{\tau}}(F + \xi(a) G) - \xi(a) a .$$
Then
$$\gamma'(a) = \beta'(\xi(a))\xi'(a) - \xi'(a)a - \xi(a) = -\xi(a) ,$$
and $\gamma(a) \leq 0$ is a concave function with strict maximum 0 at $a = \int_{R^{\tau}}G dm_F.$
In the following we assume that $G$ and $\stt$ are flow independent, which guarantees that $g(x)$ and $\tau(x)$ are $\sigma-$ independent. 
Consequently, the function $g_a = g - a \tau$ is not cohomologous to a function in $C(R: 2 \pi \Z)$,  and this yields
\be \label{eq:1.3}
\frac{d^2 \Pr(f + t g_a)}{dt^2}\Big\vert_{t = \xi(a)}=\omega(a) > 0.
\ee
 
 From now on for simplicity of the notation we will write $\Pr$ instead of $\Pr_{\sigma}$. Consider the rate function
$$J(a) =: \inf_{t \in \R} \{ \Pr(f + t(g - \tau a)\} = \Pr(f + \eta(a)(g - \tau a)),$$
where $\eta(a)$ is the unique real number such that
$$0 = \frac{d\Pr(f + t(g- \tau a))}{dt}\Bigl\vert_{t = \eta(a)}= \int_R g dm_{f+ \eta(a) (g - \tau a)} - a \int_R\tau dm_{f + \eta(a) (g - \tau a)}.$$
Notice that
$$a = \frac{\int_R \int_0^{\tau(x)}G(\pi(x, t)) dt dm_{f + \eta(a) g_a}}{\int_R \tau dm_{f + \eta(a)g_a}}= \int_{R^{\tau}}G(\pi(x, t)) dm_{F + \eta(a)(G -a)} $$
$$= \int_{R^{\tau}}G(\pi(x, t)) dm_{F +\eta(a) G} = \frac{d \Pr_{\sigma_{\tau}}(F +t G)}{dt}\Bigr \vert_{t = \eta(a)}.$$
Here we have used the fact that $F + \eta(a) G$ and $F + \eta(a) G - \eta(a) a$ have the same equilibrium state in $R^{\tau}$. 
Since $\frac{d \Pr_{\sigma_{\tau}}(F + t G)}{dt}$ is increasing, there exists an unique $\xi(a)$ such that
$$\frac{d \Pr_{\sigma_{\tau}}(F + \xi(a) G)}{dt} = a ,$$
therefore $\xi(a) = \eta(a).$ Hence $J(a) = \Pr(f + \xi(a)(g - a \tau)).$ In Section 2 we show that 
\be \label{eq:1.4}
J(a) = \gamma(a) \int_R \tau d\mu_{f + \xi(a)(g - a \tau)}.
\ee

This implies $J(a) \leq 0$ and $J(a) = 0$ if  only if $a = \int_{R^{\tau}} G dm_F$ and $\xi(a) = 0$.
We prove the following large deviation result.

\begin{thm}
Let the assumptions of Theorem $1$ be satisfied. 
Assume that $G : \Lambda \longrightarrow (0,\infty)$ is a H\"older continuous function for which
there exists a Markov family $\rr = \{ R_i\}_{i=1}^k$ for the flow $\varphi_t$ on $\mt$ such that
$G$ is constant on the stable leaves of  all "rectangular boxes"
$$B_i = \{ \varphi_t(x) : x\in R_i, 0 \leq t \leq \tau(x)\} ,$$ $i= 1, \ldots, k$. 
Assume in addition that $G$ and $\stt$ are flow independent.  Let $0 < \rho < 1$ be 
the constant in Theorem $1$ and let $\delta_n = e^{-\delta n}, \: 0 <\delta \leq -\frac{\log \rho}{2}$. Then
\begin{equation} \label{eq:1.5}
\mu\{x: g^n(x) - \tau^n(x)a \in (-\delta_n, \delta_n)\}\sim  \frac{2 \delta_n}{\sqrt{2 \pi \omega(a) n}} e^{n J(a)},\: n \to \infty.
\end{equation}
\end{thm}

A similar result for the measure of points $x \in R$ for which the difference
$$\frac{1}{n}\int_0^{\tau^n(x)} G(\pi(t, x))dt - p$$
 stays in a exponentially shrinking interval has been obtained in \cite{PeS1} under the conditions
that $G$ is a Lipschitz function on $\Lambda$ and $\frac{{\rm Lip}\:G}{\min G}< \mu_0$ with a suitable positive constant $\mu_0.$
In the present paper we improve the result in \cite{PeS1} assuming $G$ only H\"older. Moreover, here one examines a more natural difference 
$\frac{1}{\tau^n(x)}\int_0^{\tau^n(x)} G(\pi(t, x))dt - a$. This progress is essentially based on the spectral estimates for the Ruelle operator with two 
complex parameters established in Theorem 1.  A further improvement will be the analysis of the asymptotic of
$$\mu\Bigl\{x: \frac{1}{T}\int_0^T G(\pi(t, x)) dt - a \in \Bigl(-\frac{e^{\delta T}}{T}, \frac{e^{-\delta T}}{T}\Bigr)\Bigr\},\:\delta > 0,$$
as $T \to +\infty$  and this is an interesting open problem. On the other hand, the case when the interval $\Bigl(-\frac{e^{\delta T}}{T}, \frac{e^{-\delta T}}{T}\Bigr)$ 
is replaced by an interval $\Bigl(\frac{\alpha}{T}, \frac{\beta}{T}\Bigr),\: \alpha < \beta,$ has been studied in \cite{W}.
 Comparing $(\ref{eq:1.5})$  with Theorem $1$ in \cite{W}, one observes that in the case we deal with the variable tending to 
$+\infty$ by a scaling can take the form $T_n = n \int_R \tau d \mu_{\tau + \xi(a) (g- a \tau)}.$ Setting
$$\omega(a) = \frac{1}{C^2(a)}\beta''(\xi(a)) \int_R \tau d \mu_{\tau + \xi(a) (g- a \tau)},\: C(a) \neq 0,$$
one may write the leading term in $(\ref{eq:1.5})$ as
$$\frac{2 \delta_n C(a)} {\sqrt{2\pi \beta''(\xi(a))T_n}} e^{T_n \gamma(a)}$$
which modulo the constant $C(a)$ is similar to the asymptotic in \cite{W} with $T_n \to\infty$, where the rate function is precisely $\gamma(a).$

\begin{rem} The result stated in Theorem $2$ holds if we assume that $G$ is non-lattice and $g$ and $\tau$ are $\sigma-$independent. The condition $G > 0$ is 
not a restriction since we can replace $G$ by  $G + C > 0$ for some large constant $C > 0$. Then $a = \int_{R^{\tau}} G dm_F + C$, and the 
asymptotic $(\ref{eq:1.5})$ is independent on  the constant $C$. The assumption that $G$ is constant on stable leaves of rectangular boxes $B_i$
is significant, however it seems difficult to remove  when very sensitive asymptotics such as $(1.5)$ are obtained. For "standard" large deviation results,
this assumption is not necessary, since one can use Sinai's Lemma (see e.g. Proposition $1.2$ in \cite{PP}) to replace an arbitrary H\"older $G$ 
by a cohomologous function which is constant on stable leaves. In \cite{W} and \cite{PoS1}, where instead of $(-e^{-\delta n}, e^{-\delta n})$ the authors deal 
with significantly larger intervals,  however still smaller than $(-c/n, c/n)$ for a constant $c > 0$,
claims have been made that the general case of H\"older functions on two-sided shifts is easily derived from the one for one-sided shifts. 
However in both papers there are no proofs of these claims.  For sharp estimates similar to $(1.5)$, it is tempting to believe that such claims would be 
difficult to justify.
\end{rem}

\section{Preliminaries}
\renewcommand{\theequation}{\arabic{section}.\arabic{equation}}
\setcounter{equation}{0}

As in section 1, let  $\varphi_t: M \longrightarrow M$ be a $C^2$ Axiom A flow on a Riemannian manifold $M$, and let $\mt$ be a basic set 
for $\varphi_t$. The restriction of the flow on 
$\Lambda$ is a hyperbolic flow \cite{PP}. For any $x \in M$ let $W_{\epsilon}^s(x), W_{\epsilon}^u(x)$ be the local stable and
unstable manifolds through $x$,  respectively (see \cite{B2}, \cite{KH},  \cite{PP}). 
When $M$ is compact and $M$ itself is a basic set, $\varphi_t$ is called an {\it Anosov flow}.
It follows from the hyperbolicity of $\mt$  that if  $\epsilon_0 > 0$ is sufficiently small,
there exists $\ep_1 > 0$ such that if $x,y\in \mt$ and $d (x,y) < \ep_1$, 
then $W^s_{\ep_0}(x)$ and $\varphi_{[-\ep_0,\ep_0]}(W^u_{\ep_0}(y))$ intersect at exactly 
one point $[x,y ] \in \mt$  (cf. \cite{KH}). That is, there exists a unique 
$t\in [-\ep_0, \ep_0]$ such that $\varphi_t([x,y]) \in W^u_{\ep_0}(y)$. Setting $\Delta(x,y) = t$, 
defines the so called {\it temporal distance function}.

We will use the set-up and some arguments from \cite{St1} and \cite{PeS2}.  As in these papers, fix a 
{\it (pseudo) Markov family} $\rr = \{ R_i\}_{i=1}^k$ of {\it pseudo-rectangles} 
$R_i = [U_i  , S_i ] =  \{ [x,y] : x\in U_i, y\in S_i\}$. Set $R = \cup_{i=1}^{k} R_i,\: U = \cup_{i=1}^k U_i$. Consider the
{\it Poincar\'e map} $\pp: R \longrightarrow R$, defined by  $\pp(x) = \varphi_{\tau(x)}(x) \in R$, where
$\tau(x) > 0$ is the smallest positive time with $\varphi_{\tau(x)}(x) \in R$ ({\it first return time}  function). 
The  {\it shift map}  $\sigma : U   \longrightarrow U$ is given by
$\sigma  = \piU \circ \pp$, where $\piU : R \longrightarrow U$ is the {\it projection} along stable leaves. 

The hyperbolicity of the flow on $\mt$ implies the existence of constants $c_0 \in (0,1]$ and $\gamma_1 > \gamma_0 > 1$ such that
\begin{equation} \label{eq:2.1}
 c_0 \gamma_0^m\; d (u_1,u_2) \leq  d (\sigma^m(u_1), \sigma^m(u_2)) \leq \frac{\gamma_1^m}{c_0} d (u_1,u_2)
\end{equation}
whenever $\sigma^j(u_1)$ and $\sigma^j(u_2)$ belong to the same  $U_{i_j}$  for all $j = 0,1 \ldots,m$.

Define a $k \times k$ matrix $A = \{A(i, j)\}_{i, j = 1}^k$ by
$$A(i, j) = \begin{cases} 1 \:\:{\rm if}\: \pp({\rm Int}\:R_i)\cap {\rm Int}\: R_j \neq \emptyset,\\
0 \:\:{\rm otherwise}. \end{cases}$$
It is possible to construct a Markov family ${\mathcal R}$ so that $A$ is irreducible and aperiodic (see \cite{B2}).

Consider the suspension space $R^{\tau} = \{(x, t) \in R \times \R:\: 0 \leq t \leq \tau(x)\}/ \sim,$ 
where by $\sim$ we identify the points $(x, \tau(x))$ and  $(\sigma x, 0).$ 
The corresponding suspension  flow  is defined by $\sigma_t^{\tau}(x, s) = (x, s + t)$  on $R^{\tau}$ 
taking into account the identification $\sim.$
 For a H\"older continuous function $f$ on $R$, the {\it topological pressure} $\Pr(f)$ with respect to $\sigma$ is defined as
$$\Pr(f) = \sup_{m \in {\mathcal M}_{\sigma}} \big\{ h(\sigma, m) + \int f d m\big \},$$
where ${\mathcal M}_{\sigma}$ denotes the space of all $\sigma$-invariant Borel probability measures and $h(\sigma, m)$ is 
the  entropy of $\sigma$ with respect to $m$.  We say that $f$ and $g$ are {\it cohomologous} and we denote this by 
$f \sim g$ if there exists a continuous function $w$ such that $f = g + w \circ \sigma - w.$ 

The proof of (\ref{eq:1.4}) follows from the following computation:
$$\gamma(a) =\Pr_{\sigma_{\tau}}(F + \xi(a) G) - \xi(a) a = h(\sigma_{\tau}, m_{F + \xi(a) G}) + \int_{R^{\tau}}(F + \xi(a) G - \xi(a) a) dm_{F  +\xi(a) G}$$
$$ = h(\sigma_{\tau}, m_{F + \xi(a) G - \xi(a) a}) + \int_{R^{\tau}}(F + \xi(a) G - \xi(a)a) dm_{F  +\xi(a) G- \xi(a) a}$$
$$=\frac{h(\sigma, \mu_{f + \xi(a) (g - a \tau)})}{\int_R \tau d\mu_{f + \xi(a) (g- a\tau)}} + \frac{\int_R (f + \xi(a) (g - a \tau))d\mu_{f + \xi(a) (g - a \tau)}}{\int_R \tau d\mu_{f + \xi(a)( g - a\tau)}}$$
$$= \frac{\Pr(f + \xi(a)(g - a\tau))}{\int_R \tau d\mu_{f + \xi(a)( g - a\tau)}}= \frac{J(a)}{\int_R \tau d\mu_{f + \xi(a)( g - a\tau)}} .$$

\section{Proof of Theorem 2}
\renewcommand{\theequation}{\arabic{section}.\arabic{equation}}
\setcounter{equation}{0}


In this section we prove Theorem 2 exploiting the spectral estimates obtained in Theorem 1. 
We work under the assumptions of Theorem 2,  in particular, $G$ is constant on stable leaves of 
rectangular boxes $B_i$ for a certain Markov family $\rr = \{R_i\}_{i=1}^k$. Then the function
$g(x)$ depends only on $x \in U.$

We may replace $f$ by a H\"older function $\tilde{f}$ 
depending only on $x \in U$ so that with some H\"older function $z(x)$ we have
$$f(x) = \tilde{f}(x) + z(\sigma(x)) - z(x).$$
Therefore for all $t \in \R$ we have $\Pr(f + t(g - a\tau)) = \Pr(\tilde{f} + t(g - a\tau))$ and $\mu_{f} = \mu_{\tilde{f}}.$
Below we use again the notation $f$ assuming that $f(x)$ depends only on $x \in U.$

We will examine the sequence 
\be
\rho(n) = \int_{U} \chi_n(g_a^n(x)) \; d\mu\;,
\ee
where  $\chi \in C_0^{\infty}(\R: \R^+)$ is a {\it fixed  cut-off function} and
\be
\chi_n(t) = \chi(\delta_n^{-1} t) \quad , \quad x\in \R\;.
\ee

\begin{prop} Under the assumptions of Theorem $2$ we  have the asymptotic
\begin{equation}\label{eq:2.1bis}
\rho(n) \sim \frac{\delta_n}{\sqrt{2 \pi \omega(a) n}}\Bigl(\int \chi(t) dt\Bigr) e^{n J(a)},\: n \to \infty.
\end{equation}
\end{prop}

\ms
\noindent
{\it Proof.}
The Ruelle operator $\lc_{f + \xi(a)g_a}$ has a simple eigenvalue 
$$\lambda_{a} = e^{\Prf(f + \xi(a) g_a)} = e^{J(a)},$$
and so for all sufficiently small $u \in \C$ the operator $\lc_{f + (\xi(a) + i u)g_a}$ has a simple eigenvalue
$e^{\Prf(f + (\xi(a)+ \i u )g_a)}$ and the rest of the spectrum of $\lc_{f + (\xi(a) + iu)g_a}$ is contained in 
a disk of radius $\theta \lambda_{a}$ with some $0 < \theta < 1.$ Note that
$$\frac{d^2 \Pr(f + (\xi(a) + \i u))g_a)}{d^2 u}\Bigl\vert_{u = 0} =  -\omega(a) < 0.$$
Clearly for the Fourier transform $\hat{\chi}_n$ of $\chi_n$ we get
$\hat{\chi}_n(u) = \delta_n \hat{\chi}(\delta_n u).$  
Set $\omega_n(y) = e^{-\xi(a) y}\chi_n(y).$ Since $\Pr(f) = 0$, the Ruelle operator $L_f$ has a simple eigenvalues 1 and the
 adjoint operator $L_{f}^*$ satisfies 
$$L_{f}^* \mu = \mu,$$
where we denote $\mu = \mu_f$ as in Section 1.\\

Using this property and applying the Fourier transform, we have
$$\rho(n) = \frac{1}{2 \pi} \int_{-\infty}^{\infty} \Bigl( \int e^{iu g_a^n(x)}  d\mu (x)\Bigr) \hat{\chi}_n(u) du$$
$$=\frac{1}{2 \pi} \int_{-\infty}^{\infty} \Bigl( \int e^{(\xi(a) + iu) g_a^n(x)}  d\mu (x)\Bigr) \hat{\omega}_n(u) du$$
$$ = \frac{1}{2\pi} \int_{-\infty}^{\infty} \Bigl( \int L_{f + (\xi(a) + i u)g_a}^n 1 (x)  d\mu(x)\Bigr) \hat{\omega}_n(u) du$$
$$= \frac{\delta_n}{2\pi} \int_{-\infty}^{\infty} \Bigl( \int L_{f + (\xi(a) + i u)g_a}^n 1 (x) \;
d\mu(x)\Bigr) \hat{\chi}(\delta_n (u - i \xi(a))) du.$$
By Taylor expansion for small $|u|$  one gets 
$$e^{\Pr(f+ (\xi(a) + i u)g_a)} = \lambda_a\Bigr( 1 - \frac{\omega(a)}{2}u^2 + {\mathcal O}(|u|^3)\Bigr).$$
We choose $\ep_0 > 0$ sufficiently small and changing the coordinates on $(-\ep_0, \ep_0)$ by $u=\frac{\sqrt{2}v}{\sqrt{\omega(a)}}$, we write

$$I_1(n) = \frac{\delta_n}{\sqrt{2\omega(a)}\pi}\lambda_a^n \int_{-\ep_1}^{\ep_1} 
\Bigl( (1 -v^2 + iQ(v))^n(1 + {\mathcal O}(v))(\hat{\chi}(- i\delta_n \xi(a)) + {\mathcal O}(\delta_n v))\Bigr) dv  + \delta_n{\mathcal O}(\lambda_a^n \theta^n)$$ 
 with $\epsilon_1= \sqrt{\frac{\omega(a)}{2}}\epsilon_0$ and real valued function $Q(v) = {\mathcal O}(|v|^3).$ 
 The analysis of the asymptotic of this integral is given in section 4.1 in  \cite{PoS1}.
The leading term has the form
$$\frac{\delta_n}{\sqrt{2\omega(a)}\pi}\hat{\chi}(0)\lambda_a^n\int_{-\ep_1}^{\ep_1}(1 - v^2)^n dv = 
\frac{\delta_n\lambda_a^n} {\sqrt{2\omega(a)}\pi}\hat{\chi}(0)\sqrt{\frac{\pi}{n}} + \delta_n{\mathcal O}\Bigl(\frac{\lambda_a^n}{n}\Bigr)\quad ,\quad n \to \infty . $$
 Thus we deduce
\be \label{eq:2.3}
I_1(n) \sim \frac{\delta_n}{\sqrt{2 \pi \omega(a) n}} \Bigl(\int \chi(t) dt\Bigr) e^{n J(a)} + {\mathcal O}\Bigl(\frac{\delta_n e^{n J(a)}}{n}\Bigr).
\ee
Next consider the integral
$$I_2(n) = \frac{\delta_n}{2\pi}\int_{\ep_0 < |u| \leq \frac{c}{a}}\Bigl(\int L_{f + (\xi(a) + i u)g_a}^n 1(x) d\mu(x)\Bigr) 
\hat{\chi}(\delta_n(u - i\xi(a))) du$$
with $c \gg 1$ sufficiently large. Since $g_a$ is non-lattice,  for $0 < \ep_0 \leq |u| \leq \frac{c}{a}$
the operator $L_{f + (\xi(a) + iu)g_a}$ has no eigenvalues $\lambda$ with $|\lambda| = \lambda_{a}$ 
 and the spectral radius of $L_{f + (\xi(a) + i u)g_a }$ 
is strictly less than $\lambda_{a}$. Thus, there exists $\alpha= \alpha(a, c),\: 0 < \alpha < 1$, 
such that for $n \geq N(a, c)$ we have
\begin{equation} \label{eq:2.4}
\|L^n_{f + (\xi(a) + i u) g_a}\| \leq \alpha^n\lambda_{a}^n.
\end{equation}
On the other hand,
\begin{equation} \label{eq:2.5}
|\hat{\chi}(\delta_n(u - i\xi(a)))| \leq C_k \frac{e^{c_0 \delta_n |\xi(a)|}}{\delta_n^k |u|^k}, \:|u| \geq \ep_0,\: \forall k \in \N\;,
\end{equation}
with $c_0 > 0$ depending on the support of $\chi.$ Applying (\ref{eq:2.4}) and (\ref{eq:2.5}) with 
$k = 0$, for large $n$ we get
$$I_2(n) = {\mathcal O} \Bigl(\frac{\delta_n e^{n J(a)}}{n}\Bigr).$$
Now consider
$$I_3(n) = \frac{\delta_n}{2\pi} \int_{|u| > \frac{c}{a}} \Bigl( \int L_{f + (\xi(a) + i u) g_a}^n 1(x) d\mu(x)\Bigr)  \hat{\chi}(\delta_n(u - i\xi(a))) du\;.$$
We are going to use the spectral estimates established in Theorem 1 for the Ruelle operator
$$L_{f - a(\xi(a) + i u)\tau + (\xi(a) + i u)g} = L_{f + \xi(a)(g - a \tau) - i a u \tau + i  u g} .$$
 Then  $|u| \leq \frac{1}{|a|} |au|$ and for sufficiently large $|u| \geq \frac{c}{a}$ and for every $\epsilon > 0$ we are in situation to apply the spectral estimates 
\begin{eqnarray} \label{eq:2.6}
&        &\Bigl\|L_{f - \xi(a)(g - a \tau) - i a u \tau + i u g}^n 1\Bigr\|_{\infty} \nonumber\\
& \leq & C_{\epsilon} e^{n\Pr(f + \xi(a)(g - a \tau))}\rho^n |au|^{\epsilon},\:0 < \rho < 1, |au| \geq c, \: n \in \N.
\end{eqnarray}

Fix $0 < \epsilon \leq 1/2$ and apply the estimate (\ref{eq:2.5}) with $k = 2$
and (\ref{eq:2.6}) for $\epsilon$. This gives
$$|I_3(n)| \leq \delta_n \lambda_{a}^n A_{\epsilon} e^{c_0\delta_n|\xi(a)|}\frac{\rho^n |a|^{\epsilon}}{\delta_n^2} 
\int_{|u| > \frac{c}{a}} |u|^{\epsilon - 2} du = D \delta_n \lambda_{a}^n \Bigl(\frac{\rho^n}{\delta_n^2}\Bigr).$$
Recall that we have the condition
 $$0 < \delta  \leq \alpha_0 \leq -\frac{\log \rho}{2}$$
and one deduces the inequality
$$n \log \rho + 2 \delta n - \log n \leq 0 ,$$
which leads to
$$\frac{\rho^n}{\delta_n^2} \leq \frac{1}{n},\: n \geq 1\;.$$
Thus, we conclude that 
$$I_3(n) = {\mathcal O} \Bigl(\frac{\delta_n e^{n J(a)}}{n}\Bigr).$$
Consequently,
$$\rho(n) = I_1(n) + I_2(n) + I_3(n) =   \frac{\delta_n}{\sqrt{ 2 \pi \omega(a)n}}\Bigl( \int\chi(t) dt\Bigr) e^{nJ(a)}
\Bigl( 1 + {\mathcal O}(1/\sqrt{n})\Bigr)\;,n \to \infty$$
and this completes the proof of Proposition 1. 
\endofproof

\bs

To establish Theorem 2, as in \cite{PoS1}, \cite{PeS1}, we approximate the characteristic function ${\bf 1}_{[-1, 1]}$ of the interval $[-1, 1]$ by cut-off functions.

\section{Ruelle operators -- definitions and assumptions}
\renewcommand{\theequation}{\arabic{section}.\arabic{equation}}
\setcounter{equation}{0}

\def\chU{\check{U}}

Assume as in Sect. 1 that $\varphi_t: M \longrightarrow M$  is a $C^2$ weak mixing Axiom A flow 
and $\mt$ be a basic set for $\varphi_t$. Here we work under the same assumptions as these in \cite{PeS2}.
One of these is:

\medskip

\noindent
{\sc (LNIC):}  {\it There exist $z_0\in \mt$,  $\ep_0 > 0$ and $\theta_0 > 0$ such that
for any  $\ep \in (0,\ep_0]$, any $\hz\in \mt \cap W^u_{\ep}(z_0)$  and any tangent vector 
$\eta \in E^u(\hz)$ to $\mt$ at $\hz$ with  $\|\eta\| = 1$ there exist  $\tz \in \mt \cap W^u_{\ep}(\hz)$, 
$\ty_1, \ty_2 \in \mt \cap W^s_{\ep}(\tz)$ with $\ty_1 \neq \ty_2$,
$\delta = \delta(\tz,\ty_1, \ty_2) > 0$ and $\ep'= \ep'(\tz,\ty_1,\ty_2)  \in (0,\ep]$ such that
$$|\Delta( \exp^u_{z}(v), \pi_{\ty_1}(z)) -  \Delta( \exp^u_{z}(v), \pi_{\ty_2}(z))| \geq \delta\,  \|v\| $$
for all $z\in W^u_{\ep'}(\tz)\cap\mt$  and  $v\in E^u(z; \ep')$ with  $\exp^u_z(v) \in \mt$ and
$\la \frac{v}{\|v\|} , \eta_z\ra \geq \theta_0$,   where $\eta_z$ is the parallel 
translate of $\eta$ along the geodesic in $W^u_{\ep_0}(z_0)$ from $\hz$ to $z$. }

\bs

The above condition may seem complicated at a first glance, however a careful look at it shows that it is just a rather
natural non-integrability condition.

Given $x \in \mt$, $T > 0$ and  $\delta\in (0,\ep]$ set
$$B^u_T (x,\delta) = \{ y\in W^u_{\ep}(x) : d(\varphi_t(x), \varphi_t(y)) \leq \delta \: \: , \:\:  0 \leq t \leq T \} .$$

We will say that $\varphi_t$ has a {\it regular distortion along unstable manifolds} over
the basic set $\mt$  if there exists a constant $\ep_0 > 0$ with the following properties:

\ms 

(a) For any  $0 < \delta \leq   \ep \leq \ep_0$ there exists a constant $R =  R (\delta , \ep) > 0$ such that 
$$\diam( \mt \cap B^u_T(z ,\ep))   \leq R \, \diam( \mt \cap B^u_T (z , \delta))$$
for any $z \in \mt$ and any $T > 0$.

\ms

(b) For any $\ep \in (0,\ep_0]$ and any $\rho \in (0,1)$ there exists $\delta  \in (0,\ep]$
such that for  any $z\in \mt$ and any $T > 0$ we have
$\diam ( \mt \cap B^u_T(z ,\delta))   \leq \rho \; \diam( \mt \cap B^u_T (z , \ep)) .$

\bs


In this paper we work under the following {\bf  Standing Assumptions:} 

\ms

(A) $\varphi_t$ has Lipschitz local holonomy maps over $\mt$,

\ms  

(B) the local non-integrability condition (LNIC) holds for $\varphi_t$ on $\mt$,

\ms

(C) $\varphi_t$  has a regular distortion  along unstable manifolds over the basic set $\mt$.

\ms

A rather large class of examples satisfying the conditions (A) -- (C) is provided by imposing the following {\it pinching condition}:

\ms

\noindent
{\bf (P)}:  {\it There exist  constants $C > 0$ and $\beta \geq \alpha > 0$ such that for every  $x\in M$ we have
$$\frac{1}{C} \, e^{\alpha_x \,t}\, \|u\| \leq \| d\varphi_{t}(x)\cdot u\| 
\leq C\, e^{\beta_x\,t}\, \|u\| \quad, \quad  u\in E^u(x) \:\:, t > 0 $$
for some constants $\alpha_x, \beta_x > 0$ with
$\alpha \leq \alpha_x \leq \beta_x \leq \beta$ and $2\alpha_x - \beta_x \geq \alpha$ for all $x\in M$.}

\ms

We should note that  (P) holds for geodesic flows on manifolds of strictly negative sectional curvature satisfying the so called 
$\frac{1}{4}$-pinching condition.  (P) always holds when  $\dim(M) = 3$.

\ms

{ \bf Simplifying Assumptions:} $\varphi_t$ is a $C^2$ contact Anosov flow satisfying the condition (P).

\ms

By \cite{St2},  the pinching condition  (P) implies that $\varphi_t$ has  Lipschitz local holonomy maps 
and regular distortion along unstable manifolds. This and Proposition 6.1 in \cite{St2} show that:

{\bf the Simplifying Assumptions  imply the Standing Assumptions}.
 
 \ms

{\bf Throughout we work under the Standing Assumptions}. In what follows we will use arguments
similar to those in section 4 in \cite{PeS2}, however technically they will be more complicated, since the numbers
of parameters involved will increase. E.g. where we had functions $f_{at}$, $h_{at}$, etc., depending on two parameters,
now we have to deal with functions $f_{atc}$, $h_{atc}$, etc., depending on three parameters. While
some of the arguments we use here are almost the same as corresponding arguments in \cite{PeS2} (and we
omit them), there are others that require more significant modification and we do them in some detail.

Let $\rr = \{ R_i\}_{i=1}^k$ be a  {\bf fixed Markov family} for the flow $\phi_t$ on $\mt$ consisting 
of rectangles $R_i = [U_i,S_i]$ and let  $U = \cup_{i=1}^k U_i$ (see section 2). Then (2.1) hold for some
constants $c_0 \in (0,1]$ and $\gamma_1 > \gamma_0 > 1$.
Let $\hU$ be the set of those points $x\in U$ such that $\pp^m(x)$ is not a boundary point of a rectangle
for any integer $m$. In a similar way define $\hR$.

{\bf Fix a number $\alpha > 0$ and  two real-valued functions  $f$ and $g$  in $C^\alpha(\hU)$}. 
Let  $P = P_f$ be the unique real  number so that $\Pr(f- P\, \tau) = 0$, where $\Pr$ is the topological pressure
with respect to $\sigma$.
For any  $t\in \R$ with $t \geq 1$,  let $\ft$ be the {\it average of $f$ over balls in $U$ of radius} $1/t$ obtained as follows:
fix an arbitrary extension $f\in C^\alpha(V)$ (with the  same H\"older constant), where $V$ is an open 
neighborhood of $U$ in $M$, and then take the averages in question. 
Then $\ft\in C^\infty(V)$ and:

(a) $\| f-\ft\|_\infty \leq | f|_\alpha/t^\alpha$ ;
 
(b) $\Lip(\ft) \leq \Con\; \|f\|_\infty t$ ; 

(c) For any $\beta \in (0, \alpha)$ we have $| f- \ft|_{\beta} \leq 2\, |f|_\alpha/ t^{\alpha-\beta}$.

\ms

Let $G : \mt \longrightarrow \R$ be a {\bf fixed $\alpha$-H\"older function which is constant on the stable leaves} of 
all "rectangular boxes"
$B_i = \{ \varphi_t(x) : x\in R_i, 0 \leq t \leq \tau(x)\}$, $i= 1, \ldots, k$. 

Given  a large parameter $t > 0$, define $G_t$ as above, so that $G_t$ is again constant on the stable leaves of 
all rectangular boxes $B_i$  and

$(a')$ $\| G - \Gt\|_\infty \leq | G |_\alpha/t^\alpha$ ;
 
$(b')$ $\Lip(\Gt) \leq \Con\; \|G\|_\infty t$ ; 

$(c')$ For any $\beta \in (0, \alpha)$ we have $|G- \Gt|_{\beta} \leq 2\, |G|_\alpha/ t^{\alpha-\beta}$.

In particular, for some constant $C_0 > 0$ we have $\Lip(\Gt) \leq C_0 t$.

\ms

Then define $\gt : R \longrightarrow \R$ by
\be
\gt(x) = \int_0^{\tau(x)} \Gt(\pi(x,s))\, ds .
\ee
Clearly $\gt$ is $\alpha$-H\"older and constant on stable leaves, so we can regard $\gt$ as a function on $U$.
Thus, $\gt \in C^\alpha(U)$.

Let $\lambda_0 > 0$ be the {\it largest eigenvalue of  $L_{f}$}, i.e.  $\lambda_0 = e^P$,
and let $\hnu_0$ be the (unique) probability measure on $U$ with  $L_{f}^*\hnu_0 = \lambda_0\, \hnu_0$.
Fix a corresponding (positive) eigenfunction $h_{0} \in \hC^\alpha (U)$ such that $\int_U h_{0} \, d\hnu_0 = 1$. Then 
$d\nu_0 = h_0\, d\hnu_0$ defines a {\it $\sigma$-invariant probability measure} $\nu_0$ on $U$. Setting
$$f^{(0)} = f  + \ln h_{0}(u) - \ln h_{0}(\sigma(u)) - \ln \lambda_0,$$
we have $L_{f^{(0)}}^*\nu_0 =  \nu_0$, i.e. 
$\di \int_U L_{f^{(0)}} H \, d\nu_0 =  \int_U H\, d\nu_0$ for any $H \in C(U)$ and
$L_{\fo}1 = 1$.

Given real numbers $a$, $c$ and $t$ (with $|a| +  \frac{1}{|t|}$ small and $c \in I$), denote by $\lambda_{atc}$ the 
{\it largest eigenvalue} of $L_{\ft - a \tau + c \gt}$ on $\clip(U)$ and by $h_{atc}$ the corresponding (positive) 
eigenfunction such that $\int_U h_{atc}\, d\nu_{atc} = 1$, where $\hnu_{atc}$ is the unique probability measure on  
$U$ with $L_{\ft - a \, \tau + c \gt}^*\hnu_{atc} = \lambda_{atc} \, \hnu_{atc}$.  Setting 
$d\nu_{atc} = h_{atc}\, d\hnu_{atc}$ defines a {\it $\sigma$-invariant probability measure} $\nu_{atc}$ on $U$.

Given $\theta \in (0,1)$, consider the metric $d_\theta$ on $\hU$ defined by $d_\theta (x,x) = 0$ and 
$d_\theta(x,y) = \theta^m$, where $m$ is the largest integer such that $x \neq y$ belong to the same cylinder of length $m$.
Taking $\theta  \in (0,1)$ sufficiently close to $1$ and $\beta \in (0,\alpha)$ sufficiently close to $0$ we have
$(d(x,y))^\alpha \leq \Con \; d_\theta(x,y)$  and $d_\theta(x,y) \leq \Con\; (d(x,y))^\beta$ for all $x,y \in \hU$. 
In what follows we assume that $\theta$ and $\beta$ satisfy these assumptions.

By the properties of the approximations $f_t$ and $g_t$ stated above, there exists a constant $C_0 > 0$,
depending on $f$ and $\alpha$ but independent of $\beta$,  such that
\be
\| [\ft - a \tau + c g ] - f   \|_\beta \leq C_0 \, [|a| + |c| +  1/t^{\alpha-\beta}]
\ee
for all $|a|, |c| \leq 1$ and $t \geq 1$. Next, the analyticity of pressure and the 
eigenfunction projection corresponding to the {\it maximal eigenvalue} $\lambda_{atc} = e^{\Prf(\ft - a \tau + c\gt)}$ 
of the Ruelle operator $L_{\ft - a \tau + c \gt}$ on $C^\beta(U)$ (cf. e.g. Ch. 3 in \cite{PP} or Appendix 1 in \cite{La}) 
that there exists a constant $a_0 > 0$  such that, taking $C_0 > 0$ sufficiently large, we have
\be
|\Pr(\ft - a \tau + c \gt) - P| \leq C_0 \; \left(|a| + |c| +  \frac{1}{t^{\alpha-\beta}}\right) \quad , \quad
\|h_{atc} - h_0\|_\beta \leq  C_0 \; \left(|a| + |c| +  \frac{1}{t^{\alpha-\beta}}\right)\;
\ee
for  $|a|, |c| \leq a_0$ and $1/t \leq a_0$. We take $C_0 > 0$  and $a_0> 0$ so that 
$$\lambda_0/C_0 \leq \lambda_{atc} \leq C_0\, \lambda_0 ,$$
$\|\ft\|_\infty \leq C_0$ and $1/C_0 \leq h_{atc}(u) \leq C_0$ for all $u \in U$ and all $|a|, |c|, 1/t \leq a_0$.

Given real numbers $a$, $c$ and $t$ with $|a|, |c|, 1/t \leq a_0$ consider the functions
$$\fatc = \ft - a \tau + c \gt + \ln h_{atc} - \ln (h_{atc}\circ \sigma) - \ln \lambda_{atc}$$
and the operators 
$\matc = L_{\fatc} : C(U) \longrightarrow C(U).$
One checks that $\matc \; 1 = 1$.

Taking the constant $C_0 > 0$ sufficiently large, we may assume that
\be
\|\fatc - f^{(0)}\|_\beta \leq C_0\, \left[ |a| + |c| +  \frac{1}{t^{\alpha- \beta}} \right] \quad, \quad |a|, |c|, 1/t \leq a_0 .
\ee

The proof of the following lemma is given in \cite{PeS2} when  $c = 0$. In the case with three parameters
the proof is almost the same, so we omit it.

\ms

\begin{lem} Taking the constant $a_0 > 0$ sufficiently small, there exists a constant $T' > 0$ such that 
for all $a,t, c \in \R$ with $|a|, |c| \leq a_0$ and $t \geq 1/a_0$ we have $h_{atc} \in \clip(\hU)$ and
$\Lip(h_{atc}) \leq T' t$.
\end{lem}

\ms

Consequently, assuming $a_0 > 0$ is chosen sufficiently small, there exists a constant $T > 0$ 
(depending on $|f|_\alpha$ and $a_0$) such that 
\be
\|\fatc\|_\infty \leq T \quad , \quad \|\gt\|_\infty \leq T \quad , \quad\Lip(h_{atc}) \leq T\, t
\quad , \quad \Lip(\fatc) \leq T\, t
\ee
for $|a|, |c|, 1/t \leq a_0$. 
In what follows {\bf we assume that $a_0$, $C_0$, $T \geq \max \{ \, \|\tau \|_0 \, , \, \Lip(\tau_{|\hU}) \, \}$
$1 < \gamma_0 < \gamma_1$ are fixed constants}   with (2.1) and (4.2) -- (4.5).

\bs

Next,  Ruelle operators of the form $L_{f -  s \tau + z g}$, where $s = a + \i b$ and $z = c + \i w$, 
$a,b, c, w\in \R$, and $|a|, |c|\leq a_0$  for some constant $a_0 > 0$, will be studied approximating them 
by Ruelle operators of the form
$$\labtz = L_{\fatc - \i\,b\tau + z \gt} : C^{\alpha}(\hU) \longrightarrow C^\alpha (\hU) .$$
Since $\fatc - \i b \tau + z \gt$ is Lipschitz, the operators $\labtz$ preserve each of the spaces $C^{\alpha'}(\hU)$ for $0 < \alpha' \leq 1$ including the space
$\clip (\hU)$ of Lipschitz functions $h: \hU \longrightarrow \C$. For such $h$ we will denote by 
$\Lip(h)$ the Lipschitz constant of $h$.  For  $|b| \geq 1$, define the norm $\|.\|_{\lip,b}$ on $\clip (\hU)$ by 
$\| h\|_{\lip,b} = \|h\|_0 + \frac{\lip(h)}{|b|}$. Recall the norm $\|h\|_{\beta,b} = \|h\|_\infty + \frac{|h|_\beta}{|b|}$ on $C^\beta(U)$ defined in section 1. 

The main step in proving Theorem 1 is  the following.

\begin{thm} 
Under the assumptions in Theorem 1  there exist constants $0 < \rho < 1$, $a_0 > 0$, 
$b_0 \geq 1$, $A_0 > 0$ and  $C = C(B, \epsilon)> 0$ such that if $a,c\in \R$  satisfy $|a|, |c| \leq a_0$,   then
$$\|L_{\fatc - \i b \tau + (c+\i w) \gt}^m h \|_{\lip,b}  \leq C \, \rho^m \, \| h\|_{\lip,b}$$
for all $h\in \clip(\hU)$, all  integers $m \geq 1$ and all $b, w, t\in \R$ with  $|b| \geq b_0$,
$t e^{A_0t} \leq |b|$ and $|w| \leq B \, |b|$.

\end{thm}

\ms


Throughout we work under the Standing Assumptions made above and with fixed
real-valued functions $f,g \in C^\alpha(\hU)$ as in section 1, where $\alpha > 0$ is a fixed number.
Another fixed number $\beta \in (0,\alpha)$ will be used later.



Assuming that all rectangles $R_i$ are sufficiently small we have $\diam(U_i) < 1$ for all $i$.
Recall the metric $D$ on $\hU$ defined in \cite{St1}: 
$D(x,y) = \min \{ \diam(\cc) : x,y\in \cc\:, \: \cc \: \mbox{\rm a cylinder contained in }\, U_i \}$
if $x,y \in U_i$ for some $i = 1, \ldots,k$, and $D(x,y) = 1$ otherwise.  As shown in \cite{St1},  $d(x,y) \leq  D(x,y)$ 
for $x,y\in \hU_i$ for some $i$, and for any cylinder $\cc$ in $U$ the characteristic function  $\chi_{\hcc}$ of $\hcc$ on 
$\hU$ is Lipschitz with respect to $D$ and $ \Lip_D(\chi_{\hcc}) \leq 1/\diam(\cc)$.
Let $\clip_D(\hU)$ be the {\it space of all Lipschitz functions $h : \hU \longrightarrow \C$
with respect to the metric} $D$  and let $\Lip_D(h)$ be the {\it Lipschitz constant} of $h$ with respect to $D$.

Given $A > 0$, denote by $K_A(\hU)$  {\it the set of all functions } $h\in \clip_D(\hU)$  such that  $h > 0$ and
$\frac{|h(u) - h(u')|}{h(u')} \leq A\, D (u,u')$ for all $u,u' \in \hU$ that belong to the  same $\hU_i$  for some $i = 1, \ldots,k$.
For  $h\in K_A(\hU)$ we have $|\ln h(u) - \ln h(v)| \leq A\; D (u,v)$ and so
$e^{-A\; D (u,v)} \leq \frac{h(u)}{h(v)} \leq e^{A \; D (u,v)}$ for all $u, v\in \hU_i $, $i = 1, \ldots,k.$

{\bf Fix an arbitrary constant $\hgamma$ with} $ 1 < \hgamma < \gamma_0$.
The following lemma  is similar to Lemma 5.2  in \cite{PeS2}, hoverer some technical details are 
different, so we sketch its proof in the Appendix.

 
\begin{lem} Assuming $a_0 > 0$ is chosen sufficiently small, there exists a constant $A_0 > 0$ 
such that for all $a,c, t\in \R$ with $|a|, |c|\leq a_0$ and $t \geq 1$ the following hold:

\ms

(a)  {\it If $H \in K_Q(\hU)$ for some $Q > 0$, then 
$$\frac{|(\matc^m H)(u) - (\matc^m H)(u')|}{(\matc^m H)(u')} \leq
A_0 \, \left[ \frac{Q}{\hgamma^m} + e^{A_0 t} \, t \right]\, D (u,u')$$
for all $m \geq 1$ and all $u,u'\in U_i$, $i = 1, \ldots, k$.}

\ms

(b) {\it If the functions $h$ and  $H$ on $\hU$   and $Q > 0$  are such that $H > 0$ on $\hU$ and 
$$|h(v) - h(v')| \leq t\, Q\,  H(v')\, D (v,v')$$ 
for any $v,v'\in \hU_i$, $i = 1, \ldots,k$, 
then for any integer $m \geq 1$ and any $b, w \in \R$ with  $|b|,  |w|\geq 1$, for $z = c + \i w$  we have 
$$ |\labtz^m h(u) - \labtz^m h(u')| \leq  
 A_0  \left(\frac{t\, Q }{\hgamma^m} (\matc^m H)(u') + ( |b|  + e^{A_0t}t + t|w|) (\matc^m |h|)(u')\right)\, D(u,u') $$
whenever $u,u'\in \hU_i$ for some $i = 1, \ldots,k$. 
In particular,  if  $e^{A_0t} t \leq |b|$ and  $|w| \leq B |b|$ for some constant $B > 0$, then 
$$|\labtz^m h(u) - \labtz^m h(u')| \leq  
 A_1  \left(\frac{t\, Q}{\hgamma^m} (\matc^m H)(u') + t\, |b| (\matc^m |h|)(u')\right)\, D(u,u') .$$
 for some constant $A_1 > 0$, depending on $B$. }

\end{lem}

\ms


{\bf From now on we will assume that $a_0$ and $A_0$ are fixed with the properties in Lemma 2 above and
$a,b,c,w,t\in \R$ are such that $|a|, |c| \leq a_0$,  $|b|, t, |w| \geq 1$ and $|w| \leq B |b|$. As before, 
set $z = c+ \i d$.}

As in \cite{PeS2}, we need  the entire set-up and notation from section 4 in \cite{St1}, so we will now recall some of it.

Following section 4 in \cite{St1}, {\bf fix  an arbitrary point $z_0 \in \mt$ and constants  $\ep_0 > 0$ and $\theta_0 \in (0,1)$ 
 with the properties described  in} (LNIC). Assume that  $z_0 \in \Int_\mt(U_1)$, $U_1 \subset \mt \cap W^u_{\ep_0}(z_0)$ and 
$S_1 \subset \mt \cap W^s_{\ep_0}(z_0)$. {\bf Fix an arbitrary constant} $\theta_1$ such that
$0 <  \theta_0  < \theta_1 < 1 \;. $

Next, fix an arbitrary orthonormal basis $e_1, \ldots, e_{n}$ in $E^u (z_0)$ and a $C^1$
parameterization $r(s) = \exp^u_{z_0}(s)$, $s\in V'_0$, of a small neighborhood $W_0$ of $z_0$ in 
$W^u_{\ep_0} (z_0)$ such that $V'_0$ is a convex compact neighborhood of $0$ in 
$\R^{n} \approx \mbox{\rm span}(e_1, \ldots,e_n) = E^u(z_0)$. Then $r(0) = z_0$ and
$\frac{\partial}{\partial s_i} r(s)_{|s=0} = e_i$ for all $i = 1, \ldots,n$.  Set  $U'_0 = W_0\cap \mt $.
Shrinking $W_0$ (and therefore $V'_0$ as well)
if necessary, we may assume that $\overline{U'_0} \subset \Int_\mt (U_1)$ and
$\left|\left\la \frac{\partial r}{\partial s_i} (s) ,  \frac{\partial r}{\partial s_j} (s) \right\ra
- \delta_{ij}\right| $ is uniformly small for all $i, j = 1, \ldots, n$ and $s\in V'_0$, so that
$\frac{1}{2} \la \xi , \eta \ra \leq \la \; d r(s)\cdot \xi \; , \;  d r(s)\cdot \eta\;\ra \leq 2\, \la \xi ,\eta \ra$
for all $\xi , \eta \in E^u(z_0)$ and $s\in V'_0\, ,$
and
$\frac{1}{2}\, \|s-s'\| \leq  d ( r(s) ,  r(s')) \leq 2\, \|s-s'\|$ for all $s,s'\in V'_0$.

\bs

\noindent
{\bf Definitions} (\cite{St1}): (a)
For a cylinder $\cc \subset U'_0$ and a unit vector $\xi \in E^u(z_0)$
we will say that a {\it separation by a $\xi$-plane occurs} in $\cc$ if there exist $u,v\in \cc$ with 
$d(u,v) \geq \frac{1}{2}\, \diam(\cc)$ such that
$ \left\la \frac{r^{-1}(v) - r^{-1}(u)}{\| r^{-1}(v) - r^{-1}(u)\|}\;,\; \xi \right\ra  \geq \theta_1\;.$

Let  $\ss_\xi$ be the {\it family of all cylinders} $\cc$ contained in $U'_0$ such that a separation by an $\xi$-plane 
occurs in $\cc$.

\ms

(b) Given an open subset $V$ of $U'_0$  which is a finite union of open cylinders and  $\delta > 0$, let
$\cc_1, \ldots, \cc_p$ ($p = p(\delta)\geq 1$) be the family of maximal closed cylinders in $\oV$ with
$\diam(\cc_j) \leq \delta$. For any unit vector $\xi \in E^u(z_0)$ set
$M_{\xi}^{(\delta)}(V) = \cup \{ \cc_j : \cc_j \in \ss_{\xi} \:, \: 1\leq j \leq p\}\;.$

\bs

In what follows we will construct, amongst other things, a sequence of unit vectors
$\xi_1, \xi_2, \ldots, \xi_{j_0}\in E^u(z_0)$. For each $\ell = 1, \ldots,j_0$  set 
$R_\ell = \{ \eta \in \sn : \la \eta , \xi_\ell\ra  \geq \theta_0\} .$
For $t \in \R$ and $s\in E^u(z_0)$ set
$I_{\eta ,t} g(s) = \frac{g(s+t\, \eta) - g(s)}{t}$, $ t \neq 0\;$ ({\it increment} of $g$ in the direction of $\eta$). 

\ms

\begin{lem}
{\rm (\cite{St1})}
There exist integers $1 \leq n_1 \leq N_0$ and $\ell_0 \geq 1$,
a sequence of unit vectors $\eta_1, \eta_2, \ldots, \eta_{\ell_0}\in E^u(z_0)$
and a non-empty open subset $U_0$ of $U'_0$ which is a finite union of open cylinders of 
length $n_1$ such that setting $\uu = \sigma^{n_1} (U_0)$ we have:

{\rm (a)} {\it For any integer $N\geq N_0$ there exist Lipschitz maps $\vl_1, \vl_2 : U \longrightarrow U$ 
($\ell = 1,\ldots, \ell_0$)  such that $\sigma^N(\vl_i(x)) = x$  for all $x\in \uu$ and $\vl_i (\uu)$ is
a finite union of open cylinders of length $N$ ($i=1,2$; $\ell = 1,2, \ldots,\ell_0$).} 

{\rm (b)} {\it There exists a constant $\hd > 0$ such that for all  $\ell = 1, \ldots, \ell_0$, 
$s\in r^{-1}(U_0)$, $0 < |h| \leq \hd$ and $\eta \in R_\ell$ with  $s+h\, \eta \in r^{-1}(U_0\cap \mt)$ we have }
$$\left[I_{\eta,h} \left(\tau^{N}(\vl_2(\trr(\cdot ))) - \tau^{N}(\vl_1(\trr(\cdot)))\right)\right](s)  \geq \frac{\hd}{2}\,.$$

{\rm (c)} {\it We have $\overline{\vl_i (U)} \bigcap \overline{v_{i'}^{(\ell')}(U)} = \e$ whenever $(i,\ell) \neq (i',\ell')$.}

{\rm (d)} {\it For  any open cylinder $V$ in  $U_0$ there exists a constant  $\delta' = \delta'(V) > 0$  such that
$$V \subset M_{\eta_1}^{(\delta)}(V) \cup M_{\eta_2}^{(\delta)}(V) \cup \ldots  \cup M_{\eta_{\ell_0}}^{(\delta)}(V)$$ 
for all} $\delta \in (0,\delta'].$
\end{lem}


Fix  $U_0$ and $\uu$ with the properties described in Lemma 3; then $\overline{\uu} = U$.

\def\Ulo{U^{(\ell_0)}}

Set  $\di \hd = \min_{1\leq \ell\leq \ell_0} \hd_j$, $\di n_0 = \max_{1\leq \ell\leq \ell_0} m_\ell$, 
and fix an arbitrary point $\hz_0 \in U_0^{(\ell_0)}\cap \hU$. 

Fix integers $1 \leq n_1 \leq N_0$ and $\ell_0 \geq 1$, unit vectors $\eta_1, \eta_2, \ldots, \eta_{\ell_0}\in E^u(z_0)$
and a non-empty open subset $U_0$ of $W_0$ 
with the properties described in Lemma 3.
By the choice of $U_0$, $\sigma^{n_1} : U_0 \longrightarrow \uu$ is one-to-one and has an inverse
map $\psi : \uu \longrightarrow U_0$, which is Lipschitz.


Next, assume that {\bf $B > 1$, $\beta \in (0, \alpha)$ and 
$E \geq \max\left\{  4A_0\;,  \; BC_1\; , \;\frac{2 A_0\, T}{\gamma-1}\; \right\} $
are fixed constants},
where $A_0 \geq 1$ is the constant from Lemma 2 and $C_1$ is the constant from the proof of Lemma 2 in the Appendix.
 {\bf Fix an integer  $N \geq N_0$} such that
\be
\hgamma^N \geq \max \left\{ \: 6A_0 \; , \; \frac{200\, \gamma_1^{n_1}\,A_0}{c_0^2} \; , \; 
\frac{512\, \gamma_1^{n_1}\, E}{c_0\, \hd\, \rho}   \;\right\} \, .
\ee
We will also {\bf assume now that the parameter $t = t(a_0, N)> 1$ is fixed} with
\be
a_0 \leq \frac{1}{t^{\alpha-\beta}} \leq 2a_0 \quad, \quad t \leq  \frac{c_0 \hd \rho \hgamma^N}{500 E \gamma_1^{n_1} } .
\ee
(Part of this condition will be needed for the proof of Theorem 1.) Clearly the above requires to assume
that $a_0 = a_0(N)$ satisfies
$a_0^{1/(\alpha-\beta)} \leq  t .$
Some other conditions on the small parameter $a_0 = a_0(N) > 0$ will be imposed later. We will also need to choose
$$b_0 \geq t e^{A_0 t} .$$

Let the parameters $b,w\in \R$ be so that
$|w| \leq B \,   |b|$ 
and $|b|, |w| \geq b_0$.

Next,  fix maps $\vl_i : U \longrightarrow U$ ($\ell = 1, \ldots, \ell_0$, $i = 1,2$) with
the properties (a), (b), (c) and (d) in Lemma 3.  In particular, (c) gives
$\overline{v_{i}^{(\ell)}(U)} \cap \overline{v_{i'}^{(\ell')}(U)} = \e$ for all $(i,\ell) \neq (i',\ell')$.

Since $U_0$ is a finite union of open cylinders, it follows from Lemma 3(d)  that there exist a  constant
$\delta' = \delta'(U_0) > 0$  such that
$M_{\eta_1}^{(\delta)}(U_0) \cup \ldots\cup M_{\eta_{\ell_0}}^{(\delta)}(U_0) \supset U_0$
for all  $\delta\in (0, \delta']$.
{\bf Fix $\delta'$ with this property}. Set
$$\ep_1 = \min\left\{ \;\frac{1}{32 C_0 }\;,\; c_1\;,\;  \frac{1}{4E}\;,\;
\frac{1}{\hd\, \rho^{p_0+2} } \; , \; \frac{c_0r_0}{\gamma_1^{n_1}}\; , \; \frac{c_0^2(\gamma-1)}{16T\gamma_1^{n_1}}\, \right\} .$$
We will also assume that $b_0$ is chosen so that
 $\frac{\ep_1}{b_0} \leq \delta' .$

Let $\cc_m$  ($1\leq m \leq p$) be the family of {\it maximal  closed cylinders} 
contained in $\overline{U_0}$ with  $\diam(\cc_m )\leq \frac{\ep_1}{|b|}$ such that 
$U_0 \subset \cup_{j=m}^p \cc_m$ and $\overline{U_0} = \cup_{m=1}^p \cc_m$. As in \cite{St2},
\be
 \rho\, \frac{\ep_1}{|b|} \leq \diam(\cc_m ) \leq  \frac{\ep_1}{|b|} \quad, \quad 1\leq m\leq p\;.
\ee

Fix an integer $q_0 \geq 1$ such that $32 \rho^{q_0-1} < \theta_1 - \theta_0$, i.e.
$\theta_0 < \theta_1  - 32\, \rho^{q_0-1} .$
Next, let $\dd_1, \ldots, \dd_q$  be the list of all closed cylinders contained in $\overline{U_0}$ 
that are {\it subcylinders of co-length} $p_0\, q_0$ of some $\cc_m$ ($1\leq m\leq p$). 
Then  
$\overline{U_0} = \cc_1 \cup \ldots \cup \cc_p =  \dd_1 \cup \ldots \cup \dd_q .$
Moreover,
$$\rho^{p_0\, q_0+1}\cdot \frac{\ep_1}{|b|} \leq \diam(\dd_j ) \leq  \rho^{q_0}\cdot  \frac{\ep_1}{|b|} \quad, \quad 1\leq j\leq q.$$

Given $j = 1, \ldots,q $, $\ell = 1, \ldots, \ell_0$  and $i = 1,2$, set $\hdd_j = \dd_j \cap \hU$, 
$Z_j = \overline{\sigma^{n_1}(\hdd_j)}$, $\hZ_j = Z_j \cap \hU$, $\xijl = \overline{\vl_i(\hZ_j)}$,
and $\hxijl = \xijl \cap \hU$.  It then follows that $\dd_j =  \psi(Z_j)$, 
and $U = \cup_{j=1}^q Z_j$.  Moreover, $\sigma^{N-n_1}(\vl_i(x)) = \psi(x)$ for all $x \in \uu$, and
all $\xijl$ are cylinders such that $\xijl \cap {X}_{i', j'}^{(\ell')} = \e$  whenever $(i,j,\ell) \neq  (i',j', \ell')$, and
$\di \diam(\xijl) \geq \frac{c_0\, \rho^{p_0\, q_0+1}}{\gamma_1^{N}} \cdot \frac{\ep_1}{|b|}$
for all $i = 1,2$, $j = 1, \ldots,q $ and $\ell = 1, \ldots, \ell_0$. The {\it characteristic function} 
$\omega_{i,j}^{(\ell)} = \chi_{\hxijl} : \hU \longrightarrow [0,1]$ of $\hxijl$
belongs to $\clip_D(\hU)$ and $\Lip_D(\xijl) \leq 1/\diam(\xijl)$. 
Set
$$\mu_0 = \mu_0 (N) = \min \left\{\; \frac{1}{4} \; , \;  \frac{c_0  \, \rho^{p_0q_0+2}\, \ep_1}{4\, \gamma_1^N }\; , \;
\frac{1}{4\,e^{2 T N}} \,\sin^2\left(\frac{\hd\, \rho\, \ep_1}{256}\right) \; \right\} .$$

Let $J$  be a {\it subset of the set} 
$\Xi = \{\; (i,j, \ell) \; :  \; 1\leq i \leq 2\; ,\;  1\leq j\leq q\; , \;  1\leq \ell \leq \ell_0\;\}.$
Define the function  $\omega = \omega_{J} : \hU \longrightarrow [0,1]$ by
$$\di \omega  = 1- \mu_0 \,\sum_{(i, j, \ell) \in J} \eijl .$$
Clearly $\omega \in \clip_D(\hU)$ and $1-\mu_0 \leq \omega(u) \leq 1$ for any $u \in \hU$. Moreover,
$\Lip_{D} (\omega) \leq \Gamma =  \frac{ 2 \mu_0 \,\gamma_1^N}{c_0\,\rho^{p_0q_0+2}}\cdot \frac{|b|}{\ep_1} .$

Next, define {\it the contraction operator} $\nn = \nn_J(a,b,t,c) : \clip_D (\hU) \longrightarrow \clip_D (\hU)$ by 
$$\left(\nn h\right) = \matc^N (\omega_J \cdot h) .$$

Using Lemma 2 above, the proof of the following lemma is very similar to that of Lemma 5.6 in \cite{St1}
and we omit it.

\ms

\begin{lem} 
{\it Under the above conditions for $N$ and $\mu$ the following hold :}

(a) $\nn h\in K_{E|b|}(\hU)$ {\it for any} $h\in K_{E|b|}(\hU)$;

(b) {\it If $h \in \clip_D (\hU)$ and $H \in K_{E|b|}(\hU)$ are such that 
$| h (v) - h(v')|\leq E\, t\, |b| H(v')\, D (v,v')$
for any $v,v'\in U_j$, $j = 1,\ldots,k$ and $|h|\leq H$ on $\hU$ and, then for any $i = 1, \ldots,k$ and any $u,u'\in \hU_i $ we have}
$| (\labtz^N  h)(u) - (\labtz^N h)(u')| \leq E\, t\,  |b| \, (\nn H)(u')\, D (u,u') .$
\end{lem}

\ms


\noindent
{\bf Definition.}  
A subset $J$ of $\Xi$ will be called {\it dense} if for any $m = 1,\ldots, p$ 
there exists $(i,j, \ell)\in J$ such that $\dd_j \subset \cc_m$.

\ms

Denote by $\J = \J(a,b,z)$ the {\it set of all dense subsets} $J$ of $\Xi$.

\def\tgl{\tilde{\gamma}_\ell}

Although the operator $\nn$ here is different, the proof of the following lemma is very similar to that of 
Lemma 5.8 in \cite{St1} and we omit it.

\ms

\begin{lem}
Given the number $N$, there exist  $\rho_2 = \rho_2(N) \in (0,1)$ and $a_0 = a_0(N) > 0$  such that 
$$\di\int_{\hU} (\nn_J H)^2d\nu \leq \rho_{2} \,\int_{\hU} H^2 d\nu$$ 
whenever  $|a|, |c| \leq a_0$, $t \geq 1/a_0$, $J$ is dense and $H \in K_{E|b|}(\hU)$.
\end{lem}

\ms

Until the end of this section we will assume that $h, H\in \clip_D (\hU)$ are {\bf fixed functions such that }
\begin{equation}
H\in K_{E|b|}(\hU)\quad , \quad |h(u)|\leq H(u) \:\:\:\:, \:\:\; u\in \hU \;,
\ee
and
\be
|h(u) - h(u')|\leq E\, t\, |b| H(u')\, D (u,u')\:\:\:\:\:
\mbox{\rm whenever} \:\: u,u'\in \hU_i\;, \; i = 1, \ldots,k \;.
\end{equation}  

\ms

Let again $z = c + \i w$.
Define the functions $\geil : \hU  \longrightarrow \C$ ($\ell = 1, \ldots, j_0$, $i = 1,2$) by
$$\di \geol(u) = \frac{\di \left| e^{(\fatc^{N} - \i b\tau^{N} + \i w \gt^N)(\vl_1(u))} h(\vl_1(u)) +
 e^{(\fatc^{N} - \i b\tau^{N} + \i w \gt^N)(\vl_2(u))} h(\vl_2(u))\right|}{\di (1-\mu)
e^{\fatc^{N}(\vl_1(u)) }H(\vl_1(u)) + e^{\fatc^{N}(\vl_2(u))}H(\vl_2(u))} ,$$ 
$$\di \getl(u) = \frac{\di \left| e^{(\fatc^{N} - \i b\tau^{N} + \i w \gt^N)(\vl_1(u))} h(\vl_1(u)) +
 e^{(\fatc^{N} - \i b\tau^{N} + \i w \gt^N)(\vl_2(u))} h(\vl_2(u))\right|}{\di
e^{\fatc^{N}(\vl_1(u)) }H(\vl_1(u)) +  (1-\mu) e^{\fatc^{N}(\vl_2(u)) }H(\vl_2(u)) },$$ 
and set 
$$ \gl(u) = |b| \, [\tau^N(\vl_2(u)) - \tau^N(\vl_1(u))]  .$$
for all $u\in \hU $.
 
\ms

\noindent
{\bf Definitions} (\cite{St1}) 
We will say that the cylinders $\dd_j$ and $\dd_{j'}$ are {\it adjacent} if they 
are subcylinders of the same $\cc_m$ for some $m$. 
If $\dd_j$ and $\dd_{j'}$ are contained in $\cc_m$  for some $m$
and  for some  $\ell = 1, \ldots, \ell_0$ there exist  $u \in \dd_j$ and $v\in \dd_{j'}$  
such that $d(u,v) \geq \frac{1}{2}\, \diam(\cc_m)$ and 
$\left\la \frac{r^{-1}(v) - r^{-1}(u)}{\| r^{-1}(v) - r^{-1}(u)\|}\;,\; \eta_\ell \right\ra  \geq \theta_1$
we will say that $\dd_j$ and $\dd_{j'}$ are {\it $\eta_\ell$-separable in $\cc_m$}.

\bs

As a consequence of Lemma 3(b) one gets the following whose proof is almost the same
as that of Lemma 5.9 in \cite{St1}, so we omit it.

\ms

\begin{lem}
Let $j, j'\in \{ 1, 2,\ldots,q\}$ be
such that $\dd_j$ and $\dd_{j'}$ are contained in $\cc_m$ and  are $\eta_\ell$-separable 
in $\cc_m$ for some $m = 1, \ldots, p$ and $\ell = 1, \ldots, \ell_0$ . Then
$ |\gl (u) - \gl(u')| \geq A\, c_2\, \ep_1$
for all $u\in \hZ_j$ and  $u'\in \hZ_{j'}$, where $c_2 =  \frac{\hd\, \rho}{16}$.   
\end{lem}

\ms

The  following lemma is the analogue of Lemma 5.10 in \cite{St1} and represents the main step in proving Theorem 1.

\ms

\begin{lem}
Assume $|b| \geq b_0$ for some sufficiently large $b_0 > 0$, $|a|, |c| \leq a_0$, and let $|w| \leq B |b|$. 
Then  for any $j = 1, \ldots,q$  there exist $i \in \{ 1,2\}$, $j' \in \{ 1,\ldots,q\}$ and 
$\ell \in \{ 1, \ldots, \ell_0\}$  such that  $\dd_j$ and $\dd_{j'}$ 
are adjacent and $\chi^{(i)}_{\ell} (u) \leq 1$ for all  $u\in \hZ_{j'}$. 
\end{lem}

\ms

To prove this  we need the following lemma which is the analogue of  Lemma 14 in \cite{D} and its proof is very similar,
so we omit it.

\ms

\begin{lem}
 If $h$ and $H$ satisfy {\rm (4.9)-(4.10)}, then for any $j = 1, \ldots,q$, $i = 1,2$ and $\ell = 1,\ldots,  \ell_0$ we have:

(a) {\it $\frac{1}{2} \leq \frac{H(\vl_i(u'))}{H(\vl_i(u''))} \leq 2$ for all} $u', u'' \in \hZ_j$;

(b) {\it Either for all $u\in \hZ_j$ we have $|h(\vl_i(u))|\leq \frac{3}{4}H(\vl_i(u))$, or  
$|h(\vl_i(u))|\geq \frac{1}{4}H(\vl_i(u))$ for all $u\in \hZ_j$.}
\end{lem}

\ms

\noindent
{\it Sketch of proof of Lemma } 7. We use a modification of the proof of Lemma 5.10 in \cite{St1}.

Given $j = 1, \ldots, q$, let $m = 1, \ldots, p$ be such that $\dd_j \subset \cc_m$.  
As in \cite{St1} we find $j',j'' = 1, \ldots,q$ such that $\dd_{j'}, \dd_{j''} \subset \cc_m$ and 
$\dd_{j'}$ and $\dd_{j''}$ are $\eta_\ell$-separable in $\cc_m$. 

Fix $\ell$, $j'$ and $j''$ with the above properties, and set $\hZ = \hZ_j \cup \hZ_{j'}\cup \hZ_{j''}\; .$ 
If there exist $t \in \{j, j', j''\}$ and $i = 1,2$ such 
that the first alternative in Lemma 8(b) holds for $\hZ_{t}$, $\ell$  and $i$, then $\mu \leq 1/4$
implies $\chi_\ell^{(i)}(u) \leq 1$ for any $u\in \hZ_{t}$.

Assume that for every $t\in \{ j, j', j''\}$ and every $i = 1,2$ the second alternative 
in Lemma 8(b) holds for $\hZ_{t}$, $\ell$ and $i$, i.e.  $|h(\vl_i(u))|\geq \frac{1}{4}\, H(\vl_i(u))$, $u \in \hZ$.

Since $\psi(\hZ) = \hdd_j \cup \hdd_{j'} \cup \hdd_{j''} \subset \cc_m$, given $u,u'\in \hZ$ we have 
$\sigma^{N-n_1}(\vl_i(u)), \sigma^{N-n_1}(\vl_i(u'))\in \cc_m$.  Moreover, $\cc' = \vl_i (\sigma^{n_1}(\cc_m))$ is a cylinder with
$\diam(\cc') \leq \frac{\ep_1}{c_0\, \gamma_0^{N-n_1}\, |b|}$. 
Now
the estimate (6.2) in the Appendix below implies
$ |\gt^N(\vl_i(u)) - \gt^N(\vl_i(u'))| \leq  \frac{C_1 t \ep_1}{c_0\, \gamma_0^{N-n_1}\, |b|} .$
Assume for example that
$e^{c \gt^N(\vl_i(u))} |h(\vl_i(u))| \geq e^{c \gt^N(\vl_i(u))} |h(\vl_i(u'))| .$ 
Then\footnote{Using some estimates as in the proof of Lemma 2(b) in the Appendix below and 
$\|c \gt^N\|_0 \leq a_0 NT$ by (4.5).}
\begin{eqnarray*}
&      & \frac{|e^{z \gt^N(\vl_i(u))} h(\vl_i(u)) 
- e^{z \gt^N(\vl_i(u'))} h(\vl_i(u'))|}{\min\{| e^{z \gt^N(\vl_i(u))} h(\vl_i(u))| , |e^{z \gt^N(\vl_i(u'))} h(\vl_i(u'))| \}}\\
& \leq & \frac{|e^{z \gt^N(\vl_i(u))} - e^{z \gt^N(\vl_i(u'))} |}{e^{c \gt^N(\vl_i(u'))}}
+ \frac{e^{c \gt^N(\vl_i(u))} | h(\vl_i(u)) - h(\vl_i(u'))|}{ e^{c \gt^N(\vl_i(u'))} |h(\vl_i(u'))| }  \\
& \leq & \frac{|e^{z \gt^N(\vl_i(u))} - e^{z \gt^N(\vl_i(u'))} |}{e^{c \gt^N(\vl_i(u'))}}
+ \frac{ e^{c (\gt^N(\vl_i(u')) - \gt^N(\vl_i(u')))}\,E |b| H(\vl_i(u'))}{ |h(\vl_i(u'))| } D (\vl_i(u),\vl_i(u')) \\
& \leq & \frac{|e^{c \gt^N(\vl_i(u))} - e^{c \gt^N(\vl_i(u'))} |}{e^{c \gt^N(\vl_i(u'))}}
+|e^{\i w \gt^N(\vl_i(u))} - e^{\i w \gt^N(\vl_i(u'))} | +  4 E|b| {e^{2 a_0 N T}}\, \diam(\cc')\\
& \leq & (e^{C_1t} C_1 t + |w| C_1 t) \, D (\vl_i(u),\vl_i(u')) + 4 E|b| e^{2N a_0 T}\,\frac{\gamma^{n_1} \ep_1}{c_0\gamma^{N}} \\
& \leq & \frac{(B+A_0)\gamma^{n_1} \ep_1}{c_0\gamma^N} + \frac{4 E \gamma^{n_1} \ep_1}{c_0 (e^{-2a_0 T} \gamma_0)^{N}}  < \frac{\pi}{12}
\end{eqnarray*}
assuming $a_0 > 0$ is chosen sufficiently small and $N$ sufficiently large.
So, the angle between the complex numbers 
$e^{z \gt^N(\vl_i(u)} h(\vl_i(u))$ and $e^{z \gt^N(\vl_i(u')} h(\vl_i(u'))$
(regarded as vectors in $\R^2$)  is  $< \pi/6$.  In particular,  for each $i = 1,2$ we can choose a real continuous
function $\theta_i(u)$, $u \in  \hZ$, with values in $[0,\pi/6]$ and a constant $\lambda_i$ such that
$$\di e^{z \gt^N(\vl_i(u))} h(\vl_i(u)) = e^{\i(\lambda_i + \theta_i(u))} e^{c \gt^N(\vl_i(u))} |h(\vl_i(u))|$$ 
for all $u\in \hZ$. Fix an arbitrary $u_0\in \hZ$ and set $\lambda = \gamma_\ell(u_0)$. 
Replacing e.g $\lambda_2$ by $\lambda_2 +  2m\pi$ for some integer $m$, we may assume that 
$|\lambda_2 - \lambda_1 + \lambda| \leq \pi$.
Using the above, $\theta \leq 2 \sin \theta$ for $\theta \in [0,\pi/6]$, and some elementary geometry  yields
$|\theta_i(u) - \theta_i(u')|\leq 2 \sin |\theta_i(u) - \theta_i(u')| < \frac{c_2\ep_1}{8}.$

The difference between the arguments of the complex numbers
$e^{\i \,b\,\tau^N(\vl_1(u))}e^{z \gt^N(\vl_1(u)}  h(\vl_1(u))$ and 
$e^{\i \,b\, \tau^N(\vl_2(u))}e^{z \gt^N(\vl_2(u)}  h(\vl_2(u))$
is given by the function
$$\Gl(u) = [b\,\tau^N(\vl_2(u)) + \theta_2(u) + \lambda_2] -  [b\, \tau^N(\vl_1(u)) + \theta_1(u) + \lambda_1]
=  (\lambda_2-\lambda_1) + \gamma_\ell(u) + (\theta_2(u) - \theta_1(u))\;.$$
Given $u'\in \hZ_{j'}$ and $u''\in \hZ_{j''}$, since $\hdd_{j'}$ and $\hdd_{j''}$ are 
contained in $\cc_m$ and are $\eta_\ell$-separable in $\cc_m$, it follows from  Lemma 6 and the above that
\begin{eqnarray*}
|\Gl(u')- \Gl(u'')| \geq  |\gl(u') - \gl(u'')| - |\theta_1(u')-\theta_1(u'')| 
- |\theta_2(u')-\theta_2(u'')| \geq   \frac{c_2\ep_1}{2}\;.
\end{eqnarray*}
Thus,  $|\Gl(u')- \Gl(u'')|\geq \frac{c_2}{2} \epsilon_{1}$ for all $u'\in \hZ_{j'}$ and  $u''\in \hZ_{j''}$. Hence either 
$|\Gl(u')| \geq \frac{c_2}{4}\epsilon_{1}$ for all $u'\in \hZ_{j'}$ or  $|\Gl(u'')| \geq \frac{c_2}{4}\epsilon_{1}$ 
for all $u''\in \hZ_{j''}$.

Assume for example that $|\Gl(u)| \geq \frac{c_2}{4}\epsilon_{1}$ for all 
$u\in \hZ_{j'}$. Since $\hZ \subset \sigma^{n_1}(\cc_m)$, as in \cite{St1} we have
for any $u \in \hZ$ we get $|\Gamma_\ell(u)| < \frac{3\pi}{2}$.
Thus, $\frac{c_2}{4}\epsilon_{1} \leq |\Gl(u)| <  \frac{3\pi}{2}$ for all $u \in \hZ_{j'}$. Now
as in \cite{D} (see also \cite{St1}) one shows that $\chi_\ell^{(1)} (u)  \leq 1$ and $\chi_\ell^{(2)} (u)  \leq 1$ 
for all $u \in \hZ_{j'}$. 
\endofproof

\bs

Parts (a) and (b) of the following lemma can be proved in the same way as the corresponding parts of Lemma 5.3 in \cite{St1},
while part (c) follows from Lemma 4(b). 

\begin{lem}
There exist a positive integer $N$ and constants $\hrho = \hrho(N) \in (0,1)$, $a_0 = a_0(N) > 0$,
$b_0 = b_0(N) > 0$ and $E \geq 1$ such that for every $a, b,c, t \geq 1, w\in \R$ with $|a|, |c|\leq a_0$, $|b| \geq b_0$
such that $|w| \leq B |b|$, there exists a finite family $\{ \nn_J\}_{J\in \J}$ of  operators 
$$\nn_J  = \nn_J(a,b,t,c) : \clip_D (\hU) \longrightarrow \clip_D (\hU) ,$$
where $\J = \J(a,b,t,c)$, with the following properties:

(a) {\it The operators $\nn_J$ preserve the cone} $K_{E|b|} (\hU)$ ;

(b) {\it For all $H\in K_{E|b|}(\hU)$ and $J \in \J$ we have}
$\di \int_{ \hU} (\nn_J H )^2 \; d\nu_0 \leq \hrho \; \int_{\hU} H^2 \; d\nu_0$.

(c) {\it If $h, H\in \clip_D (\hU)$ are such that $H\in K_{E|b|}(\hU)$, $|h(u)|\leq H(u)$ for all   $u \in \hU$ and \\
$$| h(u) - h(u') | \leq E t |b| H(u')\, D (u,u')$$
whenever $u,u'\in \hU_i $ for some $i = 1, \ldots,k$,  
then there exists  $J\in \J$  such that $|\labtz^N h(u)|\leq (\nn_J H)(u)$ for all $u \in \hU$   and for $z= c + \i w$ we have
$$|(\labtz^N h)(u) - (\labtz^N h)(u')|\leq E\, t\, |b| (\nn_J H)(u')\, D (u,u')$$
whenever $u,u'\in \hU_i $ for some $i = 1, \ldots,k$.}
\end{lem}

\ms

\noindent
{\it Proof of Lemma } 9. 
Set $\hat{\rho} = 1- \ep_2 $.

Let $a\in \R$ and $b,w \in \R$ be such that $|a| \leq a_0$ and $|w| \leq B |b|$, $|b|, |w| \geq b_0$, and
let $J \in \J(a,b)$. Then (a) follows from  Lemma 4(a),  while  (b) follows from  Lemma 5.

To check (c), assume that $h, H \in  \clip_D (\hU)$ satisfy (4.9) and (4.10). 
Now define the subset $J$ of $\J(a,b)$ in the following way. First, include in $J$ all 
$(1,j, \ell)\in \Xi$ such that $\geol(u) \leq 1$ for all $u \in \hZ_j$. Then for any $j = 1, \ldots, q$ and 
$\ell = 1, \ldots, \ell_0$ include $(2,j,\ell)$ in $J$ if and only if $(1,j,\ell)$ has not been included in $J$ (that is, $\geol(u) > 1$ for 
some $u \in \hZ_j$) and $\getl(u) \leq 1$ for all $u \in \hZ_j$. It follows from Lemma 7  that $J$ is dense.  

Consider the operator $\nn = \nn_{J}(a,b) : \clip_D (\hU) \longrightarrow \clip_D (\hU)$. 
Then Lemma 4(b) implies 
$$| (\labtz^N \; h)(u) - (\labtz^N h)(u')| \leq E \, t\, |b| (\nn\; H)(u')\, D (u,u')$$
whenever $u,u'\in \hU_i $ for some $i = 1, \ldots,k$.  So, it remains to show that
\begin{equation}
\left| (\labtz^N h)(u)\right| \leq (\nn H)(u) \:\:\: , \:\: u\in \hU \;.\
\end{equation}

Let $u \in \hU$.  If $u \notin \hZ_j$ for any $(i,j, \ell)\in J$, then $\omega (v) = 1$ whenever
$\sigma^N v = u$ (since $v\in \xijl$ implies $u = \sigma^N v \in Z_j$). Hence
\begin{eqnarray*}
\left| (\labtz^N h)(u)\right|  =  \left| \sum_{\sigma^N v = u} e^{(\fatc^N-\i b \tau^N  + \i w \gt^N)(v)} h(v)\right|
 \leq  (\matc^N (\omega H))(u) = (\nn H)(u).
\end{eqnarray*}

Assume that $u\in \hZ_j$  for some $(i,j, \ell)\in J$. We will consider the case $ i = 1$; the case $i =2$
is similar. (Notice that by the definition of $J$, we cannot have both $(1,j, \ell)$ and  $(2,j, \ell)$ in $J$.)
Then $\geol(u) \leq 1$, and therefore
\begin{eqnarray*}
 \left| (\labtz^N h)(u)\right| 
& \leq & \left| \sum_{\sigma^N v = u, \;v\neq \vl_1(u),\vl_2(u)} e^{(\fatc^N-\i b \tau^N + \i w \gt^N)(v)} h(v) \right|\\
&        &\:\:\:       + \left| e^{(\fatc^N - \i b\tau^N + \i w \gt^N)(\vl_1(u))} h(\vl_1(u)) +  
              e^{(\fatc^N - \i b\tau^N + \i w \gt^N)(\vl_2(u))} h(\vl_2(u))\right|\nonumber\\
& \leq & \sum_{\sigma^N v = u, \;v\neq \vl_1(u),\vl_2(u)} e^{\fatc^N(v) } |h(v)|\nonumber\\
&      & + \left[(1-\mu)e^{\fatc^N  (\vl_1(u))} H(\vl_1(u)) +  e^{\fatc^N (\vl_2(u))} H(\vl_2(u))\right]\; . \nonumber
\end{eqnarray*}

Since $(1,j, \ell) \in J$ and $(2,j, \ell)\notin J$, the definition of the function $\omega$ gives
$\omega (\vl_1(u)) \geq 1-\mu$ and $\omega (\vl_2(u)) = 1$. This and (4.9) imply
\begin{eqnarray*}
\left| (\labtz^N h)(u)\right| 
&\leq  & \sum_{\sigma^N v = u, \;v\neq v_1(u),v_2(u)} e^{\fatc^N(v) } \omega (v) H(v)\\
&      &  + \left[e^{\fatc^N (v_1(u))} \omega (v_1(u)) H(v_1(u)) 
+   e^{\fatc^N  (v_2(u))} \omega (v_2(u))H(v_2(u))\right]
 =  (\nn H)(u)\;,
\end{eqnarray*}
which proves (4.11).
\endofproof



\section{Proofs of Theorems}
\renewcommand{\theequation}{\arabic{section}.\arabic{equation}}
\setcounter{equation}{0}

\noindent
{\it Proof of Theorem} 3. We use an argument from \cite{D}.

Let $B > 0$ be a constant. Let $N$, $\hrho$, $a_0$, $w_0$  and $E$ be as in Lemma 9. Given
$a,b, c,w,t\in \R$ with $|a|\leq a_0$, $|b| \geq b_0$, $|w| \leq B  |b|$, let $\{ \nn_J\}_{J\in \J}$ be a finite family
of  operators having the  properties (a), (b) and (c) in Lemma 9. 

Let $h\in \clip(U)$ be such that $\| h\|_{\lip,b} \leq 1$. Then  $|h(u)| \leq 1$ for all $u\in U$ and  $\Lip(h) \leq |b|$.  
Thus, for any $u,v\in \hU_i$. $i = 1, \ldots,k$, we have $|h(u) - h(v)| \leq |b|\, d(u,v) \leq |b|\, D(u,v)$, so
$\Lip_D(h) \leq |b|$. Set $h^{(m)} = \labtz^{m N} h$. Define the sequence of functions $\{ H^{(m)}\}$ recursively by
$H^{(0)} = 1$ and $H^{(m+1)} = \nn_{J_m} H^{(m)}$, where $J_m \in \J$ is chosen by induction so
that the conclusions of Lemma 9(c) are satisfied with $h = h^{(m)}$, $H = H^{(m)}$ and $J = J_m$. 

Since $H^{(0)}\in K_{E|b|}(U)$, it follows that $H^{(m)}\in K_{E|b|}(U)$ for all $m \geq 0$. Moreover, for
$h^{(0)} = h$ we clearly have $|h^{(0)}| \leq H^{(0)}$ and
$$| h^{(0}(u) - h^{(0)}(u') | \leq |b|\, d (u,u') \leq E t |b| H^{(0)}(u')\, D (u,u')$$
whenever $u,u'\in \hU_i $ for some $i = 1, \ldots,k$. Now Lemma 9(c) implies that
$h^{(m)}$ and $H^{(m)}$ satisfy similar conditions for all $m \geq 0$. In particular, 
$|h^{(m)}| \leq H^{(m)}$ on $\hU$ for all $m$.

Using an  induction on $m$ and property (b) in Lemma 9, we get
$$\int_{\hU} (H^{(m)})^2 \; d\nu \leq \hrho \; \int_{\hU} (H^{(m-1)})^2\; d\nu \leq \hrho^{m}.$$
Hence
$$\int_U |\labtz^{mN}h|^2 \; d\nu = \int_{\hU} |\labtz^{mN}h|^2 \; d\nu 
= \int_{\hU} |h^{(m)}|^2\; d\nu \leq \int_{\hU} (H^{(m)})^2\; d\nu \leq \hrho^m$$
for all $m \geq 1$.  This proves the theorem. 
\endofproof

\bs

\def\tfa{\tilde{f}_a}
\def\tmac{ \tilde{\mm}_{ac}}
\def\hmu{\hat{\mu}}

As in \cite{D} and \cite{St1}  we need the following lemma  whose proof is the same.

\ms

\begin{lem}
Let $\beta \in (0,\alpha)$. There exists a constants $A_1 > 0$ such that for all $a,b,c, t,w\in \R$ with 
$|a|, |c|,  1/|b|, 1/t\leq a_0$ such that $|w| \leq B  |b|$, and all positive  integers $m$  and all  $h \in C^\beta (U)$ we have
$$|\labtz^m h(u) - \labtz^m h(u')| \leq A_1\left[ \frac{|h|_\beta}{\hgamma^{m\beta}}  + |b| \, (\matc^m |h| )(u')\right]\, (d (u,u'))^\beta$$
for all $u,u'\in U_i$. 
\end{lem}

\ms

We will derive Theorem 1 from Theorem 3 and Lemma 10 above.

\bs

\noindent
{\large \it Proof of Theorem} 1.
We essentially repeat the proofs of Corollaries 2 and 3 in \cite{D} (see also the Appendix in \cite{PeS2}).

Let $\ep > 0$, $B > 0$ and $\beta \in (0,\alpha)$. Take $\rho\in (0,1)$, $a_0 > 0$, $b_0 > 0$, $A_0 > 0$ and $N$  
as in Theorem 2. We will assume that $\rho \geq \frac{1}{\gamma_0}$. Let $a, b,c,w\in \R$  be such that $|a|, |c| \leq a_0$ 
and $|b| \geq b_0$. Let $t > 0$ be such that $1/t^{\alpha - \beta} \leq a_0$. Assume that $|w| \leq B |b|$ and set $z = c + \i w$.

First, as in \cite{D} one derives from Theorem 2 and Lemma 14 (approximating functions
$h \in C^\beta(\hU)$ by Lipschitz functions) that there exist constants $C_6 > 0$ and $\rho_3 \in (0,1)$ such that
\be
\|\labtz^{n} h\|_{\beta,b} \leq C_{6} |b|^\ep \rho_{3}^n \quad , \quad n \geq 0 ,
\ee
for all $h \in C^\beta(\hU)$. 



Next, given $h \in C^\beta(\hU)$, we have
$\labtz^n(h/h_{atc}) = \frac{1}{\lambda^n_{atc} \, h_{atc} } \, L_{\ft - s\, \tau + z\, \gt} h ,$
where again $ s = a+ \i b$ and $z  = c + \i w$, so by (4.3) and (4.4) we get
\begin{eqnarray*}
\|L^n_{\ft - s \tau + z \gt}h\|_{\beta, b} 
& \leq & \lambda_{atc}^n \|h_{atc}\, \labtz^n(h/h_{atc})\|_{\beta,b}\\
& \leq & \Con \, \lambda_{0}^n (e^{3C_0a_0} \rho_3)^n |b|^\ep\, \|h/h_{atc}\|_{\beta,b} 
\leq  \Con \, \lambda_0^n\, \rho_4^n \, |b|^\ep\, \|h\|_{\beta,b}\;,
\end{eqnarray*}
where $\lambda_{atc} \leq e^{3C_0a_0} \lambda_0$ and $\rho_3 e^{3C_0a_0} = \rho_4 < 1$, provided $a_0 > 0$ is small enough.

We will now approximate $L_{f- s \tau + z g}$ by $L_{\ft - s \tau +c \gt}$
in two steps. First, the above implies
\begin{eqnarray*}
\|L^n_{f - s \tau + c g + \i w \gt} h\|_{\beta,w} 
 =     \left\|L^n_{\ft - s \tau + z\gt} \left(e^{(f^n-\ft^n) + c (g^n - \gt^n)} h\right)\right\|_{\beta,b}
 \leq  C \, \lambda_0^n\, \rho_4^n \, |b|^\ep\,  \left\|e^{(f^n-\ft^n) + c (g^n -\gt^n)} h \right\|_{\beta,b} .
\end{eqnarray*}
for some constant $C > 0$. Choosing $C$ appropriately, we have $\|f-\ft\|_0 \leq C \, a_0$ and 
$|f - \ft|_\beta \leq C /t^{\alpha-\beta} \leq C$, so  $\|f^n - \ft^n\|_0 \leq n\, \|f-\ft|_0 \leq C \, n a_0$, and similarly
$|f^n - \ft^n|_\beta \leq C\,  n a_0$. Similar estimates hold for $g^n - \gt^n$.  Thus,
$\|e^{(f^n-\ft^n) + c (g^n - \gt^n)} h \|_0 \leq e^{C \; n a_0}\|h\|_0$, and
\begin{eqnarray*}
|e^{(f^n-\ft^n)  + c (g^n -\gt^n)} h |_\beta 
& \leq &  \|e^{(f^n-\ft^n) + c (g^n -\gt^n)}\|_0\, |h|_\beta + |e^{(f^n-\ft^n)  + c(g^n -\gt^n)} |_\beta \, \|h\|_\infty\\
& \leq & e^{C \;  n a_0} |h|_\beta  + e^{C \; n a_0}\, |(f^n-\ft^n) + c (g^n -\gt^n)|_\beta\, \|h\|_\infty
 \leq  C \, n \,e^{C\;  n a_0}\,\|h\|_\beta ,
\end{eqnarray*}
replacing $C$ by a larger constant where necessary. Combining this with the previous estimate gives
$\| e^{(f^n-\ft^n) + c(g^n -\gt^n)} h \|_{\beta,b} \leq C \, n \,e^{C \; n a_0}\,\|h\|_\beta ,$
so 
$$\|L^n_{f - s \tau + cg + \i w \gt} h\|_{\beta,b} \leq
C \, \lambda_0^n\, \rho_4^n \, |b|^\ep\, n \,e^{C \; n a_0}\,\|h\|_{\beta,b} .$$
Taking $a_0 > 0$ sufficiently small, we may assume that $\rho_4\, e^{C \;  a_0} < 1$. Now take an arbitrary $\rho_5$ with
$\rho_4\, e^{C \; a_0} < \rho_5 < 1$. Then we can take the constant $C_7 > 0$
so large that $n \,\rho^n_4\, e^{C n a_0} \leq C_7 \rho_5^n$ for all  integers $n \geq 1$. This gives
$\|L^n_{f - s \tau + cg + \i w \gt} h\|_{\beta,b} \leq
C_7 \, \lambda_0^n\, \rho_5^n \, |b|^\ep\, \|h\|_{\beta,b}$ for all $n \geq 0 .$
Using the latter we can write
\begin{eqnarray*}
\|L^n_{f - s \tau + z g} h\|_{\beta,b} 
 =     \left\|L^n_{f - s \tau + c g + \i w \gt} \left(e^{\i w (g^n -\gt^n)} h\right)\right\|_{\beta,b}
 \leq  C_7 \, \lambda_0^n\, \rho_5^n \, |b|^\ep\,  \left\|e^{\i w (g^n -\gt^n)} h \right\|_{\beta,b} .
\end{eqnarray*}
We have $\|e^{\i w (g^n -\gt^n)} h \|_0 = \|h\|_0$, $|g-\gt|_\beta \leq \Con/ t^{\alpha - \beta}  \leq 1$
(if $t > 1$ is sufficiently large), so
\be
|e^{\i w (g^n -\gt^n)} h |_\beta 
 \leq   \|e^{\i w (g^n -\gt^n)}\|_0\, |h|_\beta + |e^{\i w (g^n -\gt^n)} |_\beta \, \|h\|_0
 \leq   |h|_\beta + |w|\, |g^n -\gt^n|_\beta\, \|h\|_0 .
\ee
and therefore
$\|e^{\i w(g^n -\gt^n)} h \|_{\beta,b} 
= \|e^{\i w(g^n -\gt^n)} h \|_0 + \frac{|w|}{|b|} |e^{\i w(g^n -\gt^n)} h |_\beta  \leq 2  B n \|h\|_{\beta,b} .$
This yields\\
$\|L^n_{f - s \tau + z g}h\|_{\beta,b} \leq C_8\, \lambda_0^n\,  \rho_5 ^n \, |b|^\ep\, n\, \|h\|_{\beta,b} .$
Choosing $\rho_6$ with $\rho_5 < \rho_6 < 1$ and taking the constant $C_9 > C_8$ sufficiently large, we get
$\|L^n_{f - s \tau + z g}h\|_{\beta,b}  \leq C_9 \, \lambda_0^n\,  \rho_6^n \, |b|^\ep\, \|h\|_{\beta,b}$
for all integers $n \geq 0$.
\endofproof

\section{Appendix: Proof of Lemma 2}
\renewcommand{\theequation}{\arabic{section}.\arabic{equation}}
\setcounter{equation}{0}

(a) Let  $u, u' \in \hU_i $ for some $i = 1, \ldots,k$ and let $m \geq 1$ be an integer. For any $v\in \hU$ with
$\sigma^m(v) = u$, denote by $v' = v'(v)$ the unique $v'\in \hU$ in the cylinder of length $m$ containing $v$ such
that $\sigma^m(v') = u'$.  Then
\begin{eqnarray}
|\fatc^m(v) - \fatc^m(v')|  \leq  \sum_{j=0}^{m-1} |\fatc(\sigma^j(v)) - \fatc (\sigma^j(v'))| 
 \leq  \frac{ T t}{c_0\, (\gamma_0 - 1)}\, d (u,u') \leq C_1\,t\,  D (u,u') 
\end{eqnarray}
for some constant $C_1 > 0$. Similarly, 
\be
|\gt^m(v) - \gt^m(v')| \leq C_1 \, t\, D(u,u') .
\ee
If $D(u,u') = \diam(\cc')$ for some cylinder $\cc' = C[i_{m+1}, \ldots, i_p]$,then $v,v'(v) \in \cc'' = C[i_0,i_1, \ldots,i_p]$ 
for some cylinder $\cc''$ with $\sigma^m(\cc'') = \cc'$, so 
$D(v,v') \leq \diam(\cc'')  \leq \frac{1}{c_0\, \gamma_0^m} \, \diam(\cc')  = \frac{D(u,u')}{c_0\, \gamma_0^m} .$

We have
\begin{eqnarray*}
&        &\frac{|(\matc^m H)(u) - (\matc^m H)(u')|}{\matc^m H (u')} 
 =      \frac{\di \left| \sum_{\sigma^m v = u} e^{\fatc^m(v) + c \gt^m(v)}\, H(v) -  
\sum_{\sigma^m v = u} e^{\fatc^m(v') + c \gt^m(v')}\, H(v') \right|}{\matc^m H (u')} \\
& \leq & \frac{\di \left| \sum_{\sigma^m v = u} e^{\fatc^m(v)+ c \gt^m(v)}\, (H(v) -  H(v'))\right|}{\matc^m H (u')}   +
\frac{\di \sum_{\sigma^m v = u}  \left|e^{\fatc^m(v) + c \gt^m(v)}-  e^{\fatc^m(v') + c \gt^m(v')}\right| \, H(v')}{\matc^m H (u')} \\
& \leq & \frac{\di  \sum_{\sigma^m v = u} e^{\fatc^m (v) + c \gt^m(v)}\, Q\, H(v')\, D (v,v')}{\matc^m H (u')} \\
&      &  + \frac{\di \sum_{\sigma^m v = u} \left| e^{[\fatc^m(v) + c \gt^m(v)] - [\fatc^m(v')+ c \gt^m(v')]} -1 \right|  \,
e^{\fatc^m(v') + c \gt^m(v')}\, H(v')}{\matc^m H (u')} .
\end{eqnarray*}
By (6.1) and (6.2),
$|\fatc^m(v) + c \gt^m(v)] - [\fatc^m(v')+ c \gt^m(v')| \leq 2C_1 t\, D(u,u') \leq 2C_1 t ,$
which implies
$\left| e^{[\fatc^m(v) + c \gt^m(v)] - [\fatc^m(v')+ c \gt^m(v')]} -1 \right| \leq e^{2C_1t} 2C_1 t\, D(u,u') .$
A more precise estimate follows from (4.4) and (4.5):
\begin{eqnarray*}
&        &      |\fatc^m(v) + c \gt^m(v)] - [\fatc^m(v')+ c \gt^m(v')| \nonumber\\
& \leq & |\ft^m(v) - \ft^m(v)| + |P-a| \, |\tau^m(v) - \tau^m(v')|  + |(h_{atc}(v) - h_{atc}(u)) - (h_{atc}(v') - h_{atc}(u')|\nonumber\\
&      & + a_0 |\gt^m(v) - \gt^m(v')|\nonumber\\
& \leq & 2m \|\ft - f^{(0)}\|_0 + |(f^{(0)})^m(v) - (f^{(0)})^m(v')| + \Con \, D(u,u') + 4C_0 + 2m a_0 \|\gt- g\|_0\nonumber\\
& \leq & \Con \, D(u,u') + C_2 m a_0 \leq C_2 + C_2 m\, a_0 
\end{eqnarray*}
for some constant $C_2 > 0$. Assume $a_0 > 0$ is chosen so that
$e^{C_2a_0} < \gamma_0/\hgamma .$
Then
\begin{eqnarray*}
&        &        \frac{|(\matc^m H)(u) - (\matc^m H)(u')|}{\matc^m H (u')} \\
& \leq & \frac{Q\, D (u,u')}{c_0\gamma^m}\, \frac{\di  \sum_{\sigma^m v = u}  
e^{[\fatc^m(v)+ c \gt^m(v)] - [\fatc^m(v')+ c \gt^m(v')]} e^{\fatc^m(v')+ c \gt^m(v')} \,  H(v')}{\matc^m H (u')}\\
&        & \:\:\:  + e^{2C_1 t}\, \frac{\di \sum_{\sigma^m v = u} 2C_1 t\,e^{\fatc^m(v'(v))}\, H(v'(v))}{\matc^m H (u')} \\
& \leq &  e^{C_2}\,e^{C_2m a_0 } \frac{Q\, D (u,u')}{c_0\gamma^m} + 2C_1 t e^{2C_1 t}\, D (u, u')
 \leq  A_0 \, \left[ \frac{Q}{\hgamma^m} + e^{A_0 t}\, t \right]\, D (u,u') ,
\end{eqnarray*}
for some constant $A_0 > 0$ independent of $a$, $c$, $t$, $m$ and $Q$.

\ms

(b) Let $m \geq 1$ be an integer and $u,u'\in \hU_i$ for some $i = 1, \ldots,k$. 
Using the notation $v' = v'(v)$ and the constant $C_2 > 0$  from part (a) above, where $\sigma^m v = u$ and
$\sigma^mv' = u'$, and some of the estimates from the proof of part (a), we get
\begin{eqnarray*}
&      & |\labtz^m h(u) - \labtz^m h(u')|\\
&  =   &    \di \left| \sum_{\sigma^m v = u} \left( e^{\fatc^m(v) + c \gt^m(v) - \i b \tau^m(v) + \i w \gt^m(v)}\, h(v) -
  e^{\fatc^m(v') + c \gt^m(v') - \i b \tau^m(v') + \i w \gt^m(v')}\, h(v')\right)\right|\\
& \leq & \di \left| \sum_{\sigma^m v = u} e^{\fatc^m(v) + c \gt^m(v) - \i b \tau^m(v) + \i w \gt^m(v)}\, [h(v) - h(v')]\right|\\
&      &   \di  +  \sum_{\sigma^m v = u}  \left| e^{\fatc^m(v) + c \gt^m(v)} -  e^{\fat^m(v') + c \gt^m(v')}\right|\, |h(v')|\\
&      &   \di  +   \sum_{\sigma^m v = u} \left| e^{ - \i b \tau^m(v) + \i w \gt^m(v)} -
  e^{- \i b \tau^m(v') - \i w \gt^m(v')}\right|\, e^{\fatc^m(v') + c \gt^m(v')} |h(v')|\\
& \leq & \di \sum_{\sigma^m v = u} e^{\fatc^m(v) + c \gt^m(v)}\, | h(v) - h(v')|\\
&  +   & \sum_{\sigma^m v = u} \left| e^{[\fatc^m(v) + c \gt^m(v)] 
- [\fatc^m(v') + c \gt^m(v')]} - 1\right|\, e^{\fatc^m(v') + c \gt^m(v') }\, |h(v')| \\
&      & +  \sum_{\sigma^m v = u} \left(|b|\, |\tau^m(v) - \tau^m(v')| + |w|\,  |\gt^m(v) -\gt^m(v')|\right)\, 
e^{\fatc^m(v') + c \gt^m(v')} |h(v')| .
\end{eqnarray*}

Using  the constants $C_1, C_2  > 0$ from the proof of part (a) and $e^{C_2a_0} < \gamma_0/\hgamma $ we get
\begin{eqnarray*}
        \sum_{\sigma^m v = u} e^{\fatc^m(v) + c \gt^m(v)}\, | h(v) - h(v')|
& \leq & e^{C_2}\,e^{C_2m a_0}\frac{t\, Q\, D (u,u')}{c_0\gamma_0^m}\sum_{\sigma^m v = u} e^{\fatc^m(v') + c \gt^m(v')}\, H(v')\\
& \leq & \frac{e^{C_2} t\, Q}{c_0 \hgamma^m} (\matc^m H)(u') \, D(u,u') .
\end{eqnarray*}
This, implies
\begin{eqnarray*}
|\labtz^m h(u) - \labtz^m h(u')|
& \leq & \frac{e^{C_2} t\, Q}{c_0 \hgamma^m} (\matc^m H)(u') \, D(u,u') + e^{2C_1t} 2C_1 t\, D(u,u')\, (\matc^m |h|)(u')\\
&        &  + \left( \Con \, |b| + |w|C_1 \, t \right) D(u,u')
\end{eqnarray*}
Thus, taking the constant $A_0 > 0$ sufficiently large we get
$$|(\labtz^N h)(u) - (\labtz^N  h)(u')|\\
 \leq  A_0  \left(\frac{t\, Q}{\hgamma^m} (\matc^m H)(u') + ( |b| + e^{A_0t}t + t|w|) (\matc^m |h|)(u')\right)\, D(u,u') ,$$
which proves the assertion.
\endofproof

{\footnotesize

\end{document}